\documentclass[a4paper,11pt]{article}

\usepackage[normalem]{ulem} 

\usepackage{amsmath,latexsym,amssymb,amsfonts}
\usepackage{amsmath}
\usepackage[a4paper,hmargin=1.87cm]{geometry}
\usepackage[dvips]{graphicx}
\usepackage{graphicx}
\usepackage{pifont}
 \usepackage{pstricks}
\usepackage{rotating}
\usepackage{subfigure}
\usepackage{graphics}
\usepackage{hyperref}
\usepackage{multirow}
\usepackage{bigints}
\usepackage{multicol}
\usepackage{courier}
\usepackage{pstricks}
\usepackage[latin1]{inputenc}
\usepackage{indentfirst}
\usepackage{float}
\usepackage{ccaption}
\usepackage{verbatim}
\usepackage{makeidx}
\usepackage{listings}
\usepackage{dsfont}
\usepackage{subcaption}
\usepackage[latin1]{inputenc}
\usepackage[toc,page]{appendix}
\usepackage{enumerate}

\newcommand{\proof}     {\paragraph{Proof}}
\newcommand{\carre}     {\hfill$\Box$}
\numberwithin{equation}{section}
\newtheorem{defi}{Definition}[section]
\newtheorem{lem}{Lemma}[section]
\newtheorem{theo}{Theorem}[section]

\newtheorem{prop}{Proposition}[section]
\newtheorem{rem}{Remark}

\title{Statistical inference for mean-field queueing systems}
\author{
 Ioannis Lambadaris\thanks{Department of Systems and Computer Engineering, Carleton University, Ottawa, Ontario, Canada. ioannis@sce.carleton.ca}, Ahmed Sid-Ali\thanks{School of Mathematics and Statistics,
Carleton University, Ottawa, Ontario, Canada. ahmedsidali@cunet.carleton.ca}, Wei Sun \thanks{Department of Mathematics and Statistics, Concordia University, Montreal, Canada. wei.sun@concordia.ca}, Yiqiang Q. Zhao\thanks{School of Mathematics and Statistics, Carleton University, Ottawa, Ontario, Canada. zhao@math.carleton.ca} \\
}

\date{}
\allowdisplaybreaks
\begin{document}
\maketitle

\begin{abstract}
\noindent Mean-field limits have been used now as a standard tool in approximations, including for networks with a large number of nodes. Statistical inference on mean-filed models has attracted more attention recently mainly due to the rapid emergence of data-driven systems. However, studies reported in the literature have been mainly limited to continuous models. In this paper, we initiate a study of statistical inference on discrete mean-field models (or jump processes) in terms of a well-known and extensively studied model, known as the power-of-$L$ ($L\geq 2$), or the supermarket model, to demonstrate how to deal with new challenges in discrete models. We focus on system parameter estimation based on the observations of system states at discrete time epochs over a finite period. We show that by harnessing the weak convergence results developed for the supermarket model in the literature, an asymptotic inference scheme based on an approximate least squares estimation can be obtained from the mean-field limiting equation. Also, by leveraging the law of large numbers alongside the central limit theorem, the consistency of the estimator and its asymptotic normality can be established when the number of servers and the number of observations go to infinity. Moreover, numerical results for the power-of-two model are provided to show the efficiency and accuracy of the proposed estimator.

\end{abstract}


\section{Introduction} 
The origins of mean-field theory trace back to the pioneering works of Curie \cite{Curie1895} and Weiss \cite{Weiss1907} in magnetism and phase transitions. Since then, this theory has expanded across a wide array of fields to study interacting particle systems, including statistical physics \cite{Daw83, Gart88, McKean66, McKean66(2)}, biological systems \cite{Daw2017, Mel+Ben2015}, communication networks \cite{Ben+Leboud2008, Graha2000, Gra+Mel93, Gra+Mel95}, and mathematical finance \cite{Giesecke+al2015, Kley+Klu2015}.

Moreover, the application of mean-field theory in queueing systems is traced back to the work of Dobrushin and Sukhov \cite{Dobru76} and has since proliferated due to its many benefits, see, e.g., \cite{Daw+Zhao2005, Daw+Zhao2019, Mitz96, Dobru96} for further developments. Indeed, in stochastic service systems, particularly those involving multiple parallel queues, load balancing is commonly applied to enhance performance by shortening queues, reducing wait times, and increasing throughput. This balancing mechanism effectively modifies the input-output dynamics to improve the system's quality of service. When such systems are viewed as interacting systems, mean-field theory becomes a natural framework to study their behavior. By using mean-field analysis, the performance of large systems can be evaluated by examining their limiting behavior as the system size approaches infinity. In particular, the limit often reduces when solving a deterministic system known as the mean-field limit, which corresponds to a McKean-Vlasov-type stochastic differential equation (SDE) solution. In McKean-Vlasov SDEs, the coefficients depend on both the process itself and its distribution, forming a class of non-linear SDEs. The study of such equations was initiated by McKean \cite{McKean66}, inspired by Kac's work in kinetic theory \cite{Kac56}.

Although extensive literature exists on mean-field interacting systems, research on their statistical inference has only gained attention in recent years. The pioneering work in this area is Kasonga's seminal paper \cite{Kasonga90}, which addressed parameter estimation for interacting particle systems modeled by It\^{o} SDEs through a maximum likelihood approach. After Kasonga's work, interest in the topic waned for nearly two decades before reemerging with significant contributions, as seen in \cite{AMORINO2023, belo2023, Bishwal2011, Genon+al2021, Dell+Hoff2022, SHARROCK2023} and references therein. Since then, the field has grown steadily, establishing itself as a crucial area of research. This renewed interest can be attributed to novel applications of mean-field theory and the rise of new technologies enabling access to massive datasets generated by systems of interacting agents.

To date, most statistical inference studies on mean-field models, such as those mentioned above, focus on interacting diffusion systems, with limited research on statistical inference for mean-field systems with jumps. Notable exceptions include \cite{Del+Four2016} and \cite{Liu2020}, where the authors proposed asymptotic estimation for the Bernoulli interaction parameter in a system of interacting Hawkes processes as both the number of particles and time approach infinity. In particular, as discussed in \cite{Zhao2022}, statistical inference in mean-field queueing models remains largely unexplored, despite statistical inference in queueing systems being an active area of research. The reader can consult for example \cite{Asanj+al2021} where a comprehensive survey on parameter and state estimation for queueing systems across various estimation paradigms was provided, yet it does not address statistical inference for mean-field queueing systems. To fill this gap, we propose in the current paper a statistical inference scheme for the parameters governing a specific mean-field queuing system, namely the supermarket model, also known as the power of  $L \geq 2$ choices. Thus, to the best of our knowledge, our current proposal represents a novel contribution in this area.

The supermarket model was independently introduced by Vvedenskaya et al. \cite{Dobru96} and Mitzenmacher \cite{Mitz96}. It represents a system of $N$ parallel identical queues, each served by a single server with a service rate $\nu$ and infinite buffer capacity. Tasks arrive at a rate of $N\lambda$; each task is allocated $L$ queues chosen uniformly at random among the $N$ and joins the shortest one, with ties resolved uniformly. All events in this system are independent. In particular, \cite{Dobru96} and \cite{Mitz96} studied the asymptotic behavior of the system as the number of servers becomes large, showing that the process associated with queue lengths converges to a deterministic limit represented by an infinite system of ordinary differential equations (ODEs). This model and its extensions have since become widely studied due to its theoretical and practical importance; see, e.g., \cite{Graha2000, Budhi+Fried2019, Budh+Mukh+Wu2019} and the references therein. However, as noted, the statistical inference for this model remains unexplored, which is the focus of this paper.

 We propose a statistical inference scheme to estimate the arrival and service rates in a supermarket model based on aggregate data obtained from discrete observations of a moment of the system over a finite period. To this end, we propose to exploit the ODE obtained at the mean-field limit to construct an approximate least square estimator (LSE). Then, using the law of large numbers and the central limit theorem established in the literature, we show that the proposed estimator is consistent and asymptotically normal as the number of servers and observations grows large. In addition, we test our estimator on synthetic data obtained by simulation which shows the accuracy of our approach.

 It is worthwhile to provide the following remark: An intriguing general approach to statistical inference for mean-field systems was proposed in \cite{Giesecke+al2019}. This approach leverages a  ``misspecified" or limiting model, created by approximating the system through large-systems asymptotics, incorporating the law of large numbers and central limit theorem. This enables constructing an approximate likelihood function, evaluated against the data generated by the true model. The estimator is then obtained by maximizing this approximate likelihood function. A key advantage of this method is that the approximate likelihood has a conditionally Gaussian structure, due to the central limit theorem, which allows for efficient numerical evaluation of the estimator. Although one might consider using a similar method to estimate parameters for the supermarket model, the complexity of the approximate likelihood for this model complicates the analysis of the estimator's asymptotic properties. This difficulty motivates the adoption of an alternative approach, specifically, an approximate LSE scheme.

The rest of the paper is organized as follows: First, in Section \ref{model}, we recall the supermarket model and introduce the appropriate notations. We also review some well-known asymptotic results about the model, including a new technical result in Proposition \ref{explicit-sol} that will be used in the sequel. In Section \ref{stat-inf}, we introduce the inference scheme along with our LSE and prove both the consistency and asymptotic normality of the estimator. To facilitate the reading, we put all the long proofs in the appendices. Section \ref{num-exp} provides simulations demonstrating that our estimator accurately predicts the system parameters and validates the asymptotic normality result. Finally, in Section \ref{conc}, we present conclusions and open questions, followed by the bibliography.

\section{Queuing network with selection of the shortest queue among several servers}
\label{model}
\subsection{The setting}
We start by recalling the supermarket model, first introduced separately in \cite{Mitz96} and \cite{Dobru96}. Consider a network with $N$ identical queues, each with a single server of service rate $\nu$ and an infinite buffer. Tasks arrive at rate $N\lambda$, and each task is allocated $L$ queues with uniform probability among the $N$ servers and elects to join the shortest one, ties being resolved uniformly. The $L$ selected queues may coincide. All these random events are independent. Let $X^N_i(t)$ denote the length of the $i$-th queue at time $t$ and define the empirical measure process
\begin{align*}
\varrho^N_j(t):=\frac{1}{N}\sum_{i=1}^N\mathds{1}_{\{X_i^N(t)=j\}},\quad  j=0,1,\ldots,
\end{align*}
which takes values in the space $\mathcal{P}(\mathbb{Z}_+)$ of probability measures on $\mathbb{Z}_+=\{0,1,2,\dots\}$ identified with the infinite-dimensional simplex
\begin{align*}
\mathcal{S}:=\bigg\{s=\{s_j\}_{j\in\mathbb{Z}_+}\in\mathbb{R}^{\mathbb{Z}_+}_+:\sum_{j=0}^{\infty}s_j=1\bigg\},
\end{align*}
where $\mathbb{R}_+$ is the set of all non-negative real numbers. Define the subspace $\mathcal{S}^N:=\mathcal{S}\cap\frac{\mathbb{Z}_+^{\mathbb{Z}_+}}{N}$. Thus, $\varrho^N(t)\in\mathcal{S}^N$ for all $t\geq 0$. Throughout this paper, we fix a constant  $T>0$. Let $\mathcal{D}([0,T],\mathcal{S})$ be the Skorokhod space of c\`adl\`ag functions defined on $[0,T]$ with values in $\mathcal{S}$, equipped with the usual Skorokhod topology. Let $\mathcal{C}([0,T],\mathcal{S})$ be the space of continuous functions defined on $[0,T]$ with values on $\mathcal{S}$, equipped with the uniform topology.

\subsection{Law of large numbers and central limit theorem for the empirical process}
We recall now some results describing the asymptotic behavior of the supermarket model that will be used to build the statistical inference scheme and the related analysis. For $p\ge 1$, denote by $\ell_p$ the space of $p$-th summable sequences,  i.e.,
 \begin{align*}
 \ell_p=\bigg\{x=\{x_j\}_{j\in\mathbb{Z}_+}\in\mathbb{R}^{\mathbb{Z}_+}:\sum_{j=0}^{\infty}|x_j|^p<\infty \bigg\},
\end{align*}
and denote by $\|\cdot\|_p$ the norm on it. In particular, let $\ell_2$ be the space of square summable sequences equipped with the inner product
$$\
\langle x,y \rangle=\sum_{j=0}^{\infty} x_jy_j,
$$
which makes it a Hilbert space. Moreover, define its subspace
\begin{align*}
	\tilde{\ell}_2:=\bigg\{s\in\ell_2:\sum_{j=0}^{\infty}j^2s_j^2<\infty,\,\sum_{j=0}^{\infty}s_j=0 \bigg\}.
\end{align*}
Furthermore, for any $j\in\mathbb{Z}_+$, denote by $e_j\in \ell_2$ the vector with 1 at the $j$-th coordinate and 0 elsewhere.

We first state the law of large numbers established in \cite[Theorem~3.4]{Graha2000} and reformulated in \cite[Corollary~1]{Budhi+Fried2019}.
 \begin{theo}\label{LLN-Budhi}
 Suppose that $\varrho^N(0)\rightarrow \varrho_0$ in $\mathcal{S}$ as $N\rightarrow\infty$. Then $\varrho^N\rightarrow\varrho$ in probability in $\mathcal{D}([0, T ] , \mathcal{S})$, where $\varrho$ is the unique solution in $\mathcal{C}([0, T], \mathcal{S})$ to the ODE:
 \begin{align}
 \dot{\varrho}(t)=F(\varrho(t)),\quad \varrho(0)=\varrho_0,
\label{ODE-LLN-Budhi}
 \end{align}
and, for any $x\in\ell_1$,
\begin{eqnarray}\label{F-L-op}
F(x)&=&\lambda \sum_{j=0}^{\infty}\bigg[\sum_{i=1}^L {L\choose i} x^i_{j-1}\bigg(\sum_{m=j}^{\infty}x_m \bigg)^{L-i}-\sum_{i=1}^L{L\choose i} x_j^i \bigg(\sum_{m=j+1}^{\infty}x_m \bigg)^{L-i}  \bigg]e_j\nonumber\\
&&+\nu \sum_{j=0}^{\infty}[x_{j+1}-x_j]e_j.
\end{eqnarray}
 \end{theo}

Next, we state results about the fluctuations of the empirical measure process $\varrho^N$ from its law of large number limit $\varrho$. To this end, define the process
\begin{align*}
\mathcal{Z}^N(t):=\sqrt{N}(\varrho^N (t)-\varrho (t)),\quad t\in [0,T],
\end{align*}
the operator
\begin{eqnarray}
\Phi(t)&:=&\lambda\sum_{j=0}^{\infty}(e_{j+1}-e_j)(e_{j+1}-e_j)^T\bigg(\sum_{i=1}^L {L\choose i} [\varrho_j(t)]^i\bigg(\sum_{m=j+1}^{\infty}\varrho_m(t)  \bigg)^{L-i}  \bigg)\nonumber\\
&&+\nu\sum_{j=1}^{\infty}(e_{j-1}-e_j)(e_{j-1}-e_j)^T\varrho_j(t),\label{Phi-op}
\end{eqnarray}
and the map $G:\tilde{\ell}_2\times\mathcal{S}\rightarrow\tilde{\ell}_2$ by
\begin{align}
G_j(x,s):=\frac{\partial}{\partial u}F_j(s+ux)\bigg|_{u=0}, \quad j\in\mathbb{Z}_+,x\in\tilde{\ell}_2,s\in \mathcal{S}.
\label{G-map}
\end{align}
Finally, we recall the definition of cylindrical Brownian motion which is a generalization of the scalar Brownian motion to Hilbert spaces.
  \begin{defi}
 A collection of continuous real-valued stochastic processes $\{(W_t(h))_{0\leq t\leq T}: h\in\ell_2\}$ defined on a filtered probability space $(\Omega, \mathcal{F},\mathbb{P}, \{\mathcal{F}_t\})$ is called an $\ell_2$-cylindrical Brownian motion if, for every $h\in\ell_2$, $( W_t(h))_{0\leq t\leq T}$ is an $ \{\mathcal{F}_t\}$-Brownian motion with variance $t \|h\|^2_2 $ and, for all $h,k\in\ell_2$,
 \begin{align*}
 \langle W(h), W(k)\rangle_t= t \langle h,k\rangle_2, \quad t\in [0,T].
 \end{align*}
 \end{defi}

We state now the central limit theorem introduced in \cite[Theorem 2]{Budhi+Fried2019}.
 \begin{theo}
\label{rrr3}
 Suppose that $\sup_{N\in\mathbb{N}}\sum_{j=0}^{\infty}j^2\varrho^N_j(0)<\infty$ and $\varrho^N(0)\rightarrow\varrho_0$ in $\mathcal{S}$ as $N\rightarrow\infty$. Also, suppose that $\mathcal{Z}^N(0)\rightarrow z_0$ in $\ell_2$ and that
\begin{align*}
\sup_{N\in\mathbb{N}}\sum_{j=0}^{\infty}j^2(\mathcal{Z}^N_j(0))^2<\infty.
\end{align*}
Then $\mathcal{Z}^N$ converges to $\mathcal{Z}$ in distribution in $\mathcal{D}([0, T ] , \ell_2)$ as $N\rightarrow\infty$, where $\mathcal{Z}$ is the unique weak solution to the following SDE:
\begin{align}
d \mathcal{Z} (t) = G(\mathcal{Z}(t), \varrho(t)) dt + a(t) d W(t), \quad \mathcal{Z}(0)=z_0\in \tilde{\ell}_2,
\label{CLT-Budhi-eqn}
\end{align}
$G$ is defined by $(\ref{G-map})$,  $a(t)$ is the symmetric square root of the operator $\Phi(t)$ in $(\ref{Phi-op})$, i.e., $a^2(t)=\Phi(t)$, and $W$ is an $\ell_2$-cylindrical Brownian motion.
\label{CLT-Budhi}
  \end{theo}
\begin{rem}
 The stochastic integral $\int_0^ta(s)dW(s)$ represents an $\ell_2$-valued Gausssian martingale $M(t)$ given as
\begin{align}\label{er1}
M_i(t)=\sum_{j=0}^{\infty}\int_{0}^tA_{ij}(s)dB_j(s),\quad t\in [0,T],i\in\mathbb{Z}_+,
\end{align}
with $A_{ij}(s)=\langle e_i,a(s)e_j\rangle_2$ and $\{B_j\}_{j\in\mathbb{Z}_+}$ is an independent  sequence of standard Brownian motions. The well-posedness of the SDE (\ref{CLT-Budhi-eqn}) was established in \cite[Proposition~2]{Budhi+Fried2019}.
\end{rem}

\subsection{The power of two choices model}
For the sake of simplicity, we focus in this paper on a special case of the supermarket model obtained when $L=2$. This model, known as the power of two choices, was first introduced and analyzed in \cite{Mitz96}. In this case, the operator $F$ in $(\ref{F-L-op})$ takes the following explicit form:
\begin{align}\label{rrr}
F(x)=\lambda \sum_{j=0}^{\infty}\bigg[2 x_{j-1}\sum_{m=j}^{\infty}x_m -2 x_j \sum_{m=j+1}^{\infty}x_m +x^2_{j-1}-x_j^2  \bigg]e_j+\nu \sum_{j=0}^{\infty}[x_{j+1}-x_j]e_j,\quad x\in\ell_1,
\end{align}
the map $G:\tilde{\ell}_2\times\mathcal{S}\rightarrow{\ell}_2$ in $(\ref{G-map})$ takes now the explicit form:
 \begin{align*}
 G_j(x,s)=2\lambda\sum_{m=j}^{\infty} [x_{{ j}-1}s_m+s_{{ j}-1}x_m-x_{ j}s_{m+1}-s_{ j}x_{m+1}]+\nu (x_{{ j}+1}-x_{ j}),\quad j\in\mathbb{Z}_+,x\in\tilde{\ell}_2,s\in \mathcal{S},
\end{align*}
and for any $t\in [0,T]$, the operator $\Phi(t)$ in $(\ref{Phi-op})$ becomes
  \begin{eqnarray*}
\Phi(t)&=&\lambda\sum_{j=0}^{\infty}(e_{j+1}-e_j)(e_{j+1}-e_j)^T\bigg( 2 \varrho_j(t)\sum_{m=j+1}^{\infty}\varrho_m(t) +[\varrho_j(t)]^2\bigg)\\
&& +\nu\sum_{j=1}^{\infty}(e_{j-1}-e_j)(e_{j-1}-e_j)^T\varrho_j(t).
\end{eqnarray*}

The following result shows that the solution to  $(\ref{CLT-Budhi-eqn})$ is a Gaussian process. Although the proof is given for the power of two choices model, it can be adapted to cover the case with general $L\geq 2$.

\begin{prop}\label{explicit-sol}
Suppose that L = 2. Then, the solution $\mathcal{Z}(t)$ to the SDE $(\ref{CLT-Budhi-eqn})$ is a Gaussian process.
\end{prop}
	\proof See Appendix \ref{explicit-sol-proof}. \carre


\section{Statistical inference of the supermarket model}
\label{stat-inf}
Suppose the service and arrival rates, \(\nu\) and \(\lambda\), that govern the system are unknown. Our goal is to estimate these parameters using observations collected over a specific time interval \([0, T]\). The complexity of the system's dynamics makes brute-force Monte Carlo estimation computationally intensive, particularly as the number of nodes, \(N\), increases. Therefore, we propose developing an estimator that utilizes the weak convergence results outlined in Section \ref{model}, specifically the law of large numbers presented in Theorem \ref{LLN-Budhi} and the central limit theorem in Theorem \ref{CLT-Budhi}.

In particular, we construct an approximate LSE based on the ODE given in \((\ref{ODE-LLN-Budhi})\). We subsequently demonstrate that this estimator is consistent and asymptotically normal as both the system size and the number of observations tend to infinity. For convenience, we denote the vector of unknown parameters as \(\theta = (\lambda, \nu)\).

\subsection{The data}
Our objective is to estimate the parameter vector \(\theta\) that governs the dynamics of the power-of-two model based on observations collected over a finite time interval. Specifically, we assume that observations are not available for every server in the network; instead, they are gathered as an aggregate measure of the system. Collecting individual data for each server can be prohibitively costly in practice, particularly for large networks, which justifies our approach to data collection.

In this context, we assume that the available data for inference includes observations of the empirical measure of the system, \(\varrho^N(t)\), over the finite time interval \([0, T]\) at \(m\) discrete points, defined as \(t_k = \frac{kT}{m}\):
\[
\varrho^N_j(t_k) := \frac{1}{N}\sum_{i=1}^N\mathds{1}_{\{X_i^N(t_k) = j\}}, \quad j \in \mathbb{Z}_+, \quad k = 1, 2, \dots, m.
\]
The observed data is then represented as
\[
D^{N,m} := \{\varrho^N(t_k): 1 \leq k \leq m\}.
\]
Thus, this dataset reflects a realization of the system governed by the true parameters, say \(\theta^* = (\lambda^*, \nu^*)\).


\subsection{Least square estimator (LSE)}
Recall from Theorem \ref{LLN-Budhi} that the empirical measure $\varrho^N$ converges in probability towards $\varrho$ the unique solution to the ODE $(\ref{ODE-LLN-Budhi})$ as the number of servers $N\rightarrow\infty$. We propose to take advantage of this result to build an approximate LSE for the parameters $\theta^*=(\lambda^*,\nu^*)$ based on the dataset $D^{N,m}$.

\paragraph{The least square function:}
Let us first introduce the following functions defined for all $j\in\mathbb{Z}_+$ and $x\in\ell_1$,:
\begin{align}
	U_j(x):=2 x_{j-1}\sum_{i=j}^{\infty}x_i -2 x_j \sum_{i=j+1}^{\infty}x_i +x^2_{j-1}-x_j^2 ,
\label{U-func}
\end{align}
and
\begin{eqnarray}\label{Oct21b}
	V_j(x):=x_{j+1}-x_j.
\end{eqnarray}
Therefore, by $(\ref{rrr})$, one can observe that
\begin{eqnarray}\label{Oct21b}
	F(x)=\lambda \sum_{j=0}^{\infty}U_j(x)e_j+\nu \sum_{j=0}^{\infty}V_j(x)e_j.
\end{eqnarray}
Moreover, let $\varrho(t)$ be the unique solution to the ODE (\ref{ODE-LLN-Budhi}). To highlight the dependence on the parameter $\theta=(\lambda,\nu)$, we will write $\varrho(t,\theta)$ in the sequel. Furthermore, let us introduce the quadratic function defined for any $\lambda,\nu\in\mathbb{R}$ by
\begin{align*}
{\cal G}(\lambda,\nu):=\sum_{j=0}^{\infty}\left[{\varrho}_j(T,\theta^*)-{\varrho}_j(0,\theta^*)-\lambda\int_0^T U_j(\varrho(s,\theta^*))ds-\nu \int_0^TV_j(\varrho(s,\theta^*))ds\right]^2.
\end{align*}
 Notice that, by (\ref{ODE-LLN-Budhi}), one immediately observes that ${\cal G}(\lambda^*,\nu^*)=0$. Thus, the function ${\cal G}(\lambda,\nu)$ attains its minimum at $\theta^*=(\lambda^*,\nu^*)$, namely
\begin{equation}
	\begin{split}
	0&=\left.\frac{\partial {\cal G}}{\partial\lambda}\right|_{(\lambda,\nu)=(\lambda^*,\nu^*)}\\
&=-2\sum_{j=0}^{\infty}\left[{\varrho}_j(T,\theta^*)-{\varrho}_j(0,\theta^*)-\lambda^*\int_0^T U_j(\varrho(s,\theta^*))ds-\nu^* \int_0^TV_j(\varrho(s,\theta^*))ds\right]\int_0^T U_j(\varrho(s,\theta^*))ds,
\label{deriv1}
\end{split}
\end{equation}
and
\begin{equation}
	\begin{split}
	0&=\left.\frac{\partial {\cal G}}{\partial\nu}\right|_{(\lambda,\nu)=(\lambda^*,\nu^*)}\\
&=-2\sum_{j=0}^{\infty}\left[{\varrho}_j(T,\theta^*)-{\varrho}_j(0,\theta^*)-\lambda^*\int_0^T U_j(\varrho(s,\theta^*))ds-\nu^* \int_0^TV_j(\varrho(s,\theta^*))ds\right]\int_0^T V_j(\varrho(s,\theta^*))ds.
\label{deriv2}
\end{split}
\end{equation}
Solving the equations $(\ref{deriv1})$ and $(\ref{deriv2})$ leads to
\begin{equation}
	\begin{split}
	\begin{pmatrix}
		\lambda^*\\
		\nu^*
	\end{pmatrix}=\begin{pmatrix}
		a_{11}&a_{12}\\
		a_{21}&a_{22}
	\end{pmatrix}^{-1}\begin{pmatrix}
		b_1\\
		b_2
	\end{pmatrix}=\frac{1}{a_{11}a_{22}-(a_{12})^2}\begin{pmatrix}
		a_{22}&-a_{12}\\
		-a_{21}&a_{11}
	\end{pmatrix}\begin{pmatrix}
		b_1\\
		b_2
	\end{pmatrix},
	\label{lambda-nu-star}
\end{split}
\end{equation}
where
$$
a_{11}:=\sum\limits_{j=0}^{\infty}\left[\int_0^T U_j(\varrho(s,\theta^*))ds\right]^2,
$$
$$
a_{12}=a_{21}:=\sum\limits_{j=0}^{\infty}\int_0^T U_j(\varrho(s,\theta^*))ds\int_0^T V_j(\varrho(s,\theta^*))ds,
$$
$$
a_{22}:=\sum\limits_{j=0}^{\infty}\left[\int_0^T V_j(\varrho(s,\theta^*))ds\right]^2,
$$
$$
b_1:=\sum\limits_{j=0}^{\infty}\left[{\varrho}_j(T,\theta^*)-{\varrho}_j(0,\theta^*)\right]\int_0^T U_j(\varrho(s,\theta^*))ds,
$$
and
$$
b_2:=\sum\limits_{j=0}^{\infty}\left[{\varrho}_j(T,\theta^*)-{\varrho}_j(0,\theta^*)\right]\int_0^T V_j(\varrho(s,\theta^*))ds.
$$
Note that the right-hand side in $(\ref{lambda-nu-star})$ is well posed only if $a_{11}a_{22}\neq (a_{12})^2$. Nevertheless, a simple application of H\"older's inequality tells us that $a_{11}a_{22}\ge(a_{12})^2$. Therefore, one needs to investigate the conditions under which the inequality
\begin{equation}\label{equal-Hold}
a_{11}a_{22}>(a_{12})^2
\end{equation}
holds. Below we give as examples two sufficient conditions that ensure the validity of (\ref{equal-Hold}).

\begin{lem} \label{lem:3.1}
	Suppose that one of the following conditions holds:
	\begin{align}
		\int_0^T\varrho_1(s)ds> \int_0^T\varrho_0(s)ds>0,
		\label{cond-ineq-posed1}
		\end{align}
		or
		\begin{align}
			\int_0^T\varrho_0(s)ds> \int_0^T\varrho_1(s)ds> \int_0^T\varrho_2(s)ds>0.
					\label{cond-ineq-posed2}
		\end{align}
	Then, the inequality $(\ref{equal-Hold})$ holds.
\end{lem}
\proof See Appendix \ref{lem3-1-proof}. \carre

\paragraph{The approximate LSE:}
Recall that the data $D^{N,m}$ are observed on the prelimiting finite $N$-system. Hence, using $(\ref{lambda-nu-star})$ we construct the following approximate LSE defined for any $m\in\mathbb{N}$ and $N\geq 2$ as
\begin{eqnarray}\label{Nov1b}
	\begin{pmatrix}
		\lambda^{N,m}\\
		\nu^{N,m}
	\end{pmatrix}:=\frac{1}{a^{N,m}_{11}a^{N,m}_{22}-(a^{N,m}_{12})^2}\begin{pmatrix}
		a^{N,m}_{22}&-a^{N,m}_{12}\\
		-a^{N,m}_{21}&a^{N,m}_{11}
	\end{pmatrix}\begin{pmatrix}
		b^{N,m}_1\\
		b^{N,m}_2
	\end{pmatrix},
	\label{appr-est}
\end{eqnarray}
where
$$
a^{N,m}_{11}:=\sum\limits_{j=0}^{\infty}\left[\frac{T}{m}\sum\limits_{k=1}^mU_j(\varrho^N(t_k,\theta^*))\right]^2,
$$
$$
a^{N,m}_{12}=a^{N,m}_{21}:=\sum\limits_{j=0}^{\infty}\left[\frac{T}{m}\sum\limits_{k=1}^mU_j(\varrho^N(t_k,\theta^*))\right]\left[\frac{T}{m}\sum\limits_{k=1}^mV_j(\varrho^N(t_k,\theta^*))\right],
$$
$$
a^{N,m}_{22}:=\sum\limits_{j=0}^{\infty}\left[\frac{T}{m}\sum\limits_{k=1}^mV_j(\varrho^N(t_k,\theta^*))\right]^2,
$$
$$
b^{N,m}_1:=\sum\limits_{j=0}^{\infty}\left[{\varrho}^N_j(T,\theta^*)-{\varrho}^N_j(0,\theta^*)\right]\left[\frac{T}{m}\sum\limits_{k=1}^mU_j(\varrho^N(t_k,\theta^*))\right],
$$
and
$$
b^{N,m}_2:=\sum\limits_{j=0}^{\infty}\left[{\varrho}^N_j(T,\theta^*)-{\varrho}^N_j(0,\theta^*)\right]\left[\frac{T}{m}\sum\limits_{k=1}^mV_j(\varrho^N(t_k,\theta^*))\right],
$$
with $t_k=\frac{kT}{m}$ for $k=1,2,\dots, m$. We will show next that the estimator in $(\ref{appr-est})$ is consistent and approximately normal.

\subsection{Consistency of the LSE}
We show in this section the consistency of the approximate estimator given in ($\ref{Nov1b}$). To this end, let us first prove the following technical lemma:
	\begin{lem}\label{lem1} Let $\{\mu_n,n\in \mathbb{Z}_+\}$ be a sequence of probability measures on $\mathbb{Z}_+$. Then, the following statements are equivalent:
	\begin{enumerate}[(i)]
	\item $\mu_n$ converges weakly to $\mu_0$ as $n\rightarrow\infty$.
    \item $\lim_{n\rightarrow\infty}\mu_n(j)=\mu_0(j)$ for all $j\in\mathbb{Z}_+$.
    \item $\lim_{n\rightarrow\infty}\sum_{j=0}^{\infty}|\mu_n(j)-\mu_0(j)|=0$.	
	\end{enumerate}		
	\end{lem}
	\noindent {\bf Proof.} Obviously, (iii)$\Rightarrow$(i)$\Rightarrow$(ii). (ii)$\Rightarrow$(iii) is complete by the observation that, for any $m\in\mathbb{N}$,
	\begin{eqnarray*}
		\sum_{j=0}^{\infty}\left|\mu_n(j)-\mu_0(j)\right|&\le&\sum_{j=0}^m\left|\mu_n(j)-\mu_0(j)\right|+\sum_{j=m+1}^{\infty}\mu_n(j)+\sum_{j=m+1}^{\infty}\mu_0(j)\\
		&=&\sum_{j=0}^m\left|\mu_n(j)-\mu_0(j)\right|+\sum_{j=0}^m[\mu_0(j)-\mu_n(j)]+2\sum_{j=m+1}^{\infty}\mu_0(j)\\
		&\le&2\sum_{j=0}^m\left|\mu_n(j)-\mu_0(j)\right|+2\sum_{j=m+1}^{\infty}\mu_0(j).
	\end{eqnarray*}
\carre

We are now ready to prove the consistency of the estimator in (\ref{Nov1b}).
\begin{theo}\label{cons-theo} Suppose that $\varrho^N(0)\rightarrow \varrho_0$ in $\mathcal{S}$ as $N\rightarrow\infty$ and (\ref{equal-Hold}) holds. Then,
	\begin{eqnarray*}
		\begin{pmatrix}
			\lambda^{N,m}\\
			\nu^{N,m}
		\end{pmatrix}\longrightarrow
		\begin{pmatrix}
			\lambda^*\\
			\nu^*
		\end{pmatrix}\ \ {\rm in\ probability}\ {\rm as} \ N\rightarrow\infty\ {\rm and }\ m\rightarrow\infty.
	\end{eqnarray*}
	 \end{theo}
\proof See Appendix \ref{theo-3-1-proof}. \carre

\subsection{Asymptotic normality of the LSE}
We now analyze the asymptotic distribution of the approximate LSE given in $(\ref{Nov1b})$. In particular, one aims to prove the asymptotic normality of
\begin{eqnarray}
	&&\sqrt{N}\bigg( \begin{pmatrix}
		\lambda^{N,m}\\
		\nu^{N,m}
	\end{pmatrix}- \begin{pmatrix}
		\lambda^*\\
		\nu^*
	\end{pmatrix}\bigg),
	\label{est-error}
\end{eqnarray}
as the network size $N$ and the number of observations $m$ go to infinity. To do so, we will prove that (\ref{est-error}) converges towards a linear combination of the process $\mathcal{Z}$ which is itself a Gaussian process by Proposition \ref{explicit-sol}. 

To ease of notations, define the following quantities:
\begin{equation*}
		\begin{split}
	\mathcal{I}&=\sum\limits_{j=0}^{\infty}\Bigg\{ a_{22}\Bigg(\left[{\varrho}_j(T,\theta^*)-{\varrho}_j(0,\theta^*)\right]\Bigg[\sum\limits_{i=0}^{\infty}\int_0^T \partial_iU_j(\varrho(s))\mathcal{Z}_i(s,\theta^*)ds\Bigg]\\
&\quad\quad\quad\quad \quad\quad\quad\quad   +\left[\mathcal{Z}_j(T,\theta^*)-\mathcal{Z}_j(0,\theta^*)\right]\int_0^T U_j(\varrho(s))ds\Bigg)\\
	&\quad\quad + 2b_1\int_0^T V_j(\varrho(s))ds\Bigg[\sum\limits_{i=0}^{\infty}\int_0^T \partial_iV_j(s,\theta^*)\mathcal{Z}_i(s,\theta^*)ds\Bigg]\\
&\quad\quad -b_2\Bigg(\int_0^T V_j(\varrho(s))ds\Bigg[\sum\limits_{i=0}^{\infty}\int_0^T \partial_iU_j(\varrho(s))\mathcal{Z}_i(s,\theta^*)ds\Bigg]\\
&\quad\quad\quad\quad \quad +\int_0^T U_j(\varrho(s))ds\Bigg[\sum\limits_{i=0}^{\infty}\int_0^T \partial_iV_j(\varrho(s))\mathcal{Z}_i(s,\theta^*)ds\Bigg]\Bigg)\\
&\quad\quad - a_{12}\Bigg(\left[{\varrho}_j(T,\theta^*)-{\varrho}_j(0,\theta^*)\right]\Bigg[\sum\limits_{i=0}^{\infty}\int_0^T \partial_iV_j(\varrho(s))\mathcal{Z}_i(s,\theta^*)ds\Bigg]\\
&\quad\quad\quad\quad\quad  +\left[\mathcal{Z}_j(T,\theta^*)-\mathcal{Z}_j(0,\theta^*)\right]\int_0^T V_j(\varrho(s))ds\Bigg)\Bigg\},
\label{I}
\end{split}
\end{equation*}
	\begin{equation*}
		\begin{split}
	\mathcal{J}&=\sum\limits_{j=0}^{\infty}\Bigg\{ a_{11}\Bigg(\left[{\varrho}_j(T,\theta^*)-{\varrho}_j(0,\theta^*)\right]\Bigg[\sum\limits_{i=0}^{\infty}\int_0^T \partial_iV_j(\varrho(s))\mathcal{Z}_i(s,\theta^*)ds\Bigg]\\
&\quad\quad\quad\quad \quad\quad +\left[\mathcal{Z}_j(T,\theta^*)-\mathcal{Z}_j(0,\theta^*)\right]\int_0^T V_j(\varrho(s))ds\Bigg)\\
	&\quad\quad + 2b_2\int_0^T U_j(\varrho(s))ds\Bigg[\sum\limits_{i=0}^{\infty}\int_0^T \partial_iU_j(s,\theta^*)\mathcal{Z}_i(s,\theta^*)ds\Bigg]\\
&\quad\quad -b_1\Bigg(\int_0^T V_j(\varrho(s))ds\Bigg[\sum\limits_{i=0}^{\infty}\int_0^T \partial_iU_j(\varrho(s))\mathcal{Z}_i(s,\theta^*)ds\Bigg]\\
&\quad\quad\quad\quad \quad +\int_0^T U_j(\varrho(s))ds\Bigg[\sum\limits_{i=0}^{\infty}\int_0^T \partial_iV_j(\varrho(s))\mathcal{Z}_i(s,\theta^*)ds\Bigg]\Bigg)\\
&\quad\quad - a_{12}\Bigg(\left[{\varrho}_j(T,\theta^*)-{\varrho}_j(0,\theta^*)\right]\Bigg[\sum\limits_{i=0}^{\infty}\int_0^T \partial_iU_j(\varrho(s))\mathcal{Z}_i(s,\theta^*)ds\Bigg]\\
&\quad\quad\quad\quad \quad +\left[\mathcal{Z}_j(T,\theta^*)-\mathcal{Z}_j(0,\theta^*)\right]\int_0^T U_j(\varrho(s))ds\Bigg)\Bigg\},
\label{IP}
\end{split}
\end{equation*}
and
\begin{equation*}
		\begin{split}
	\mathcal{K}&= 2\sum\limits_{j=0}^{\infty}\Bigg\{-a_{11}\int_0^T V_j(\varrho(s))ds\Bigg[\sum\limits_{i=0}^{\infty}\int_0^T \partial_iV_j(\varrho(s))\mathcal{Z}_i(s,\theta^*)ds\Bigg]\\
&\quad\quad\quad -a_{22}\int_0^T U_j(\varrho(s))ds\Bigg[\sum\limits_{i=0}^{\infty}\int_0^T \partial_iU_j(\varrho(s))\mathcal{Z}_i(s,\theta^*)ds\Bigg]\\
	&\quad\quad\quad + a_{12}\Bigg(\int_0^T V_j(\varrho(s))ds\Bigg[\sum\limits_{i=0}^{\infty}\int_0^T \partial_iU_j(\varrho(s))\mathcal{Z}_i(s,\theta^*)ds\Bigg]\\
&\quad\quad\quad\quad\quad\quad +\int_0^T U_j(\varrho(s))ds\Bigg[\sum\limits_{i=0}^{\infty}\int_0^T \partial_iV_j(\varrho(s))\mathcal{Z}_i(s,\theta^*)ds\Bigg]\Bigg)\Bigg\}.
\label{IP2}
\end{split}
\end{equation*}
One can observe that $\mathcal{I}$, $\mathcal{J}$, $\mathcal{K}$ are linear combinations of the Gaussian process $\mathcal{Z}(t)$ solution to the SDE $(\ref{CLT-Budhi-eqn})$ (see Proposition \ref{explicit-sol}). 

The following theorem,  states the asymptotic normality of the approximate LSE (\ref{est-error}) under the assumption that $\frac{m}{\sqrt{N}}\rightarrow\infty$.
	
	\begin{theo}\label{thm9} Suppose  that (\ref{equal-Hold}) holds. Then, under the assumption of Theorem \ref{rrr3}, as $N,m,\frac{m}{\sqrt{N}}\rightarrow\infty$,
		\begin{align*}
			\sqrt{N}\bigg( \begin{pmatrix}
				\lambda^{N,m}\\
				\nu^{N,m}
			\end{pmatrix}- \begin{pmatrix}
				\lambda^*\\
				\nu^*
			\end{pmatrix}\bigg)&\xrightarrow{\text{d}}\frac{1}{[a_{11}a_{22}-(a_{12})^2]^2}\begin{pmatrix}
				[a_{11}a_{22}-(a_{12})^2]\mathcal{I}+(a_{22}b_1-a_{12}b_2)\mathcal{K}\\
				\\
				[a_{11}a_{22}-(a_{12})^2]\mathcal{J}+(a_{11}b_2-a_{12}b_1)\mathcal{K}
			\end{pmatrix}.
		\end{align*}
			\end{theo}
			\proof See Appendix \ref{theo-2-proof}. \carre

\begin{rem}
	By Proposition \ref{explicit-sol}, we know that the limiting distribution given by Theorem \ref{thm9} is normal. Moreover, explicit (but tedious) expressions of the mean and the covariance matrix of the limiting normal distribution can be obtained by virtue of (\ref{er1}) and (\ref{5}).
\end{rem}

\section{Numerical experiments}
\label{num-exp}
\begin{figure}
    \centering
    \subfigure[]{\includegraphics[width=0.4\textwidth]{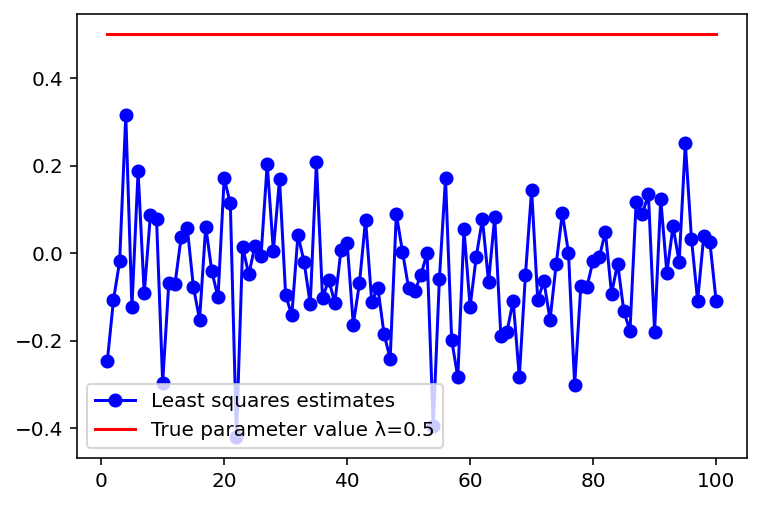}}
    \subfigure[]{\includegraphics[width=0.4\textwidth]{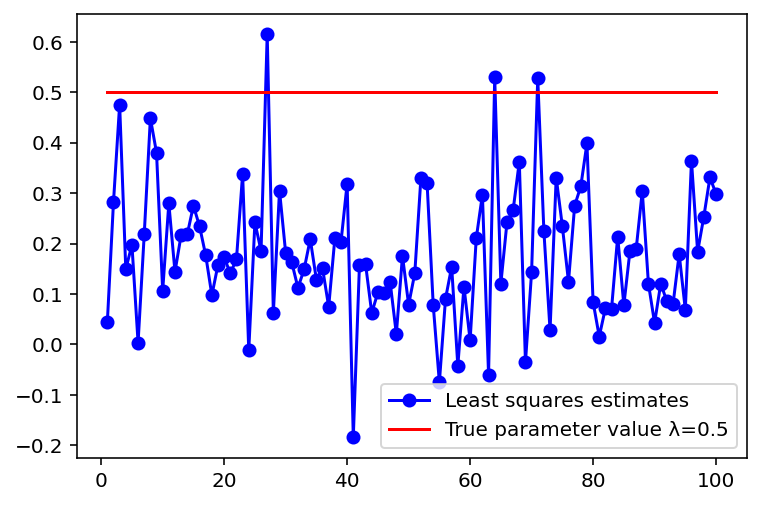}}
    \subfigure[]{\includegraphics[width=0.4\textwidth]{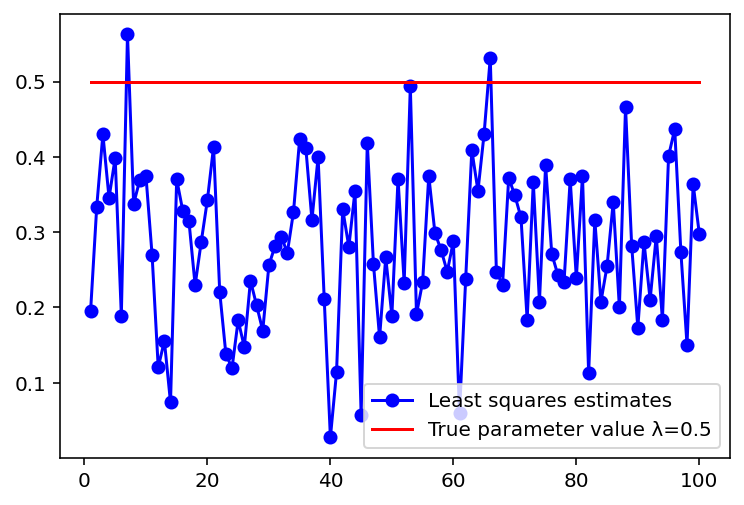}}
    \subfigure[]{\includegraphics[width=0.4\textwidth]{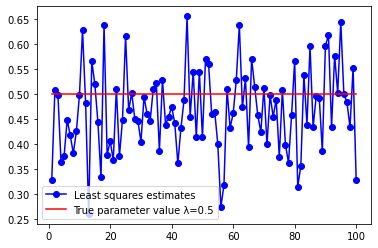}}
    \subfigure[]{\includegraphics[width=0.4\textwidth]{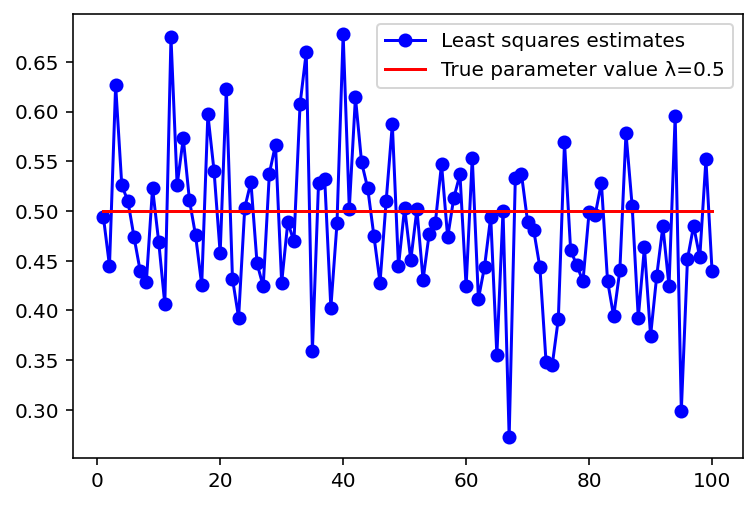}}
    \caption{Values of $100$ Least squares estimates of the parameter $\lambda=0.5$ with different network sizes $N$ and number of observations $m$: (a) N=100, m=1000; (b) N=500, m=10000; (c) N=1000, m=10000; (d) N=2000, m=20000; (e) N=3000, m=30000}
    \label{est_lamb}
\end{figure}

\begin{figure}
    \centering
    \subfigure[]{\includegraphics[width=0.4\textwidth]{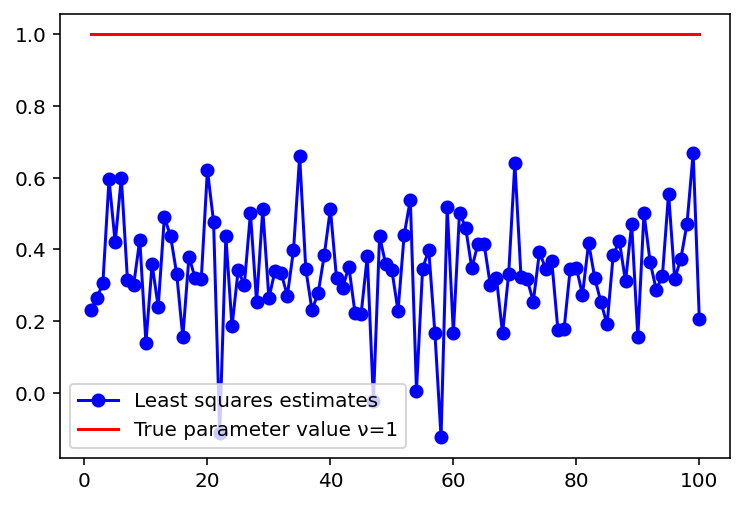}}
    \subfigure[]{\includegraphics[width=0.4\textwidth]{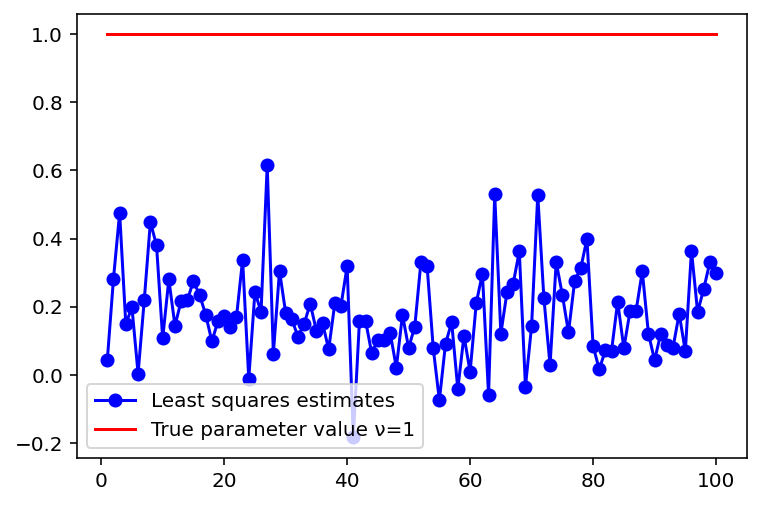}}
    \subfigure[]{\includegraphics[width=0.4\textwidth]{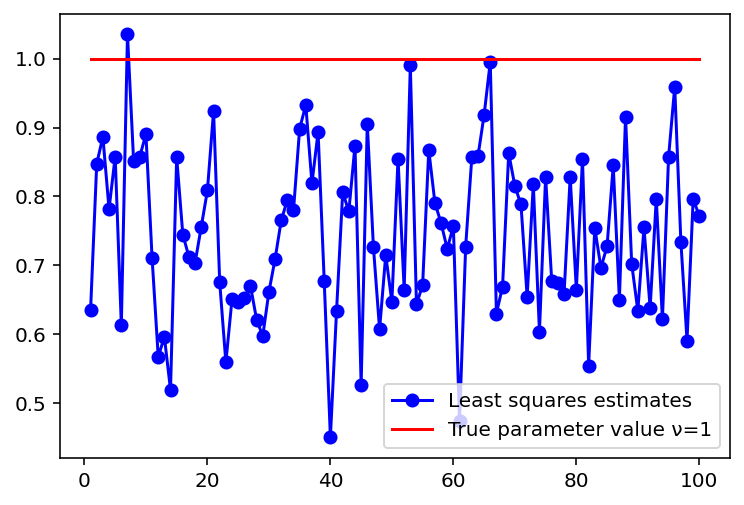}}
    \subfigure[]{\includegraphics[width=0.4\textwidth]{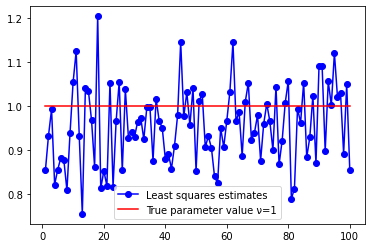}}
    \subfigure[]{\includegraphics[width=0.4\textwidth]{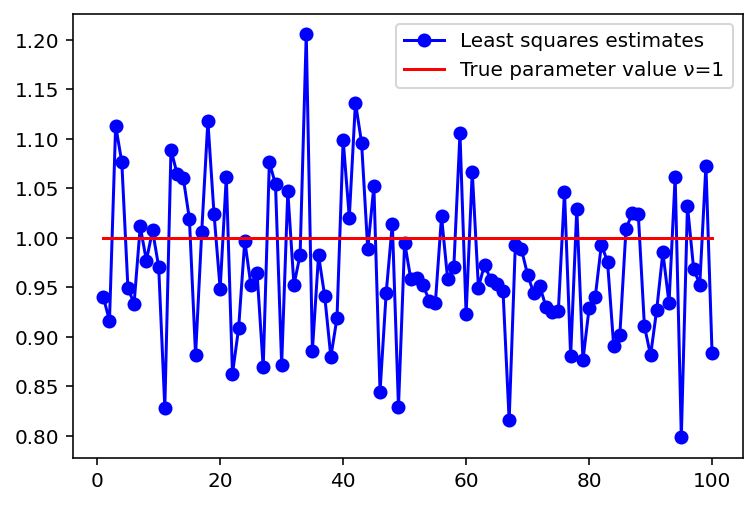}}
    \caption{Values of $100$ Least squares estimates of the parameter $\nu=1$ with different network sizes $N$ and number of observations $m$: (a) N=100, m=1000; (b) N=500, m=10000; (c) N=1000, m=10000; (d) N=2000, m=20000; (e) N=3000, m=30000}
    \label{est_nu}
\end{figure}

In this section, we evaluate the consistency and asymptotic normality of the approximate LSE $\hat{\theta}=(\lambda^{N,m},\nu^{N,m})$ defined in $(\ref{appr-est})$ using simulated data. Specifically, we aim to validate the assertions made in Theorem $\ref{cons-theo}$ and Theorem $\ref{thm9}$. To this end, we propose to generate the datasets $D^{N,m}=\{\varrho^N(t_k): 1 \leq k \leq m\}$ by simulating the power of two choices model, with true parameters $(\theta=(\lambda,\nu)=(0.5,1))$, for various network sizes $N\in\mathbb{N}$ and different numbers of observations $m\in\mathbb{N}$. For each combination of $(N,m)$, we simulate a total of $100$ samples which are then utilized to estimate the arrival and service rates $\theta=(\lambda,\nu)$ using the approximate LSE $\hat{\theta}=(\lambda^{N,m},\nu^{N,m})$. The estimated values are plotted in Figure \ref{est_lamb} and Figure \ref{est_nu}.

\subsection{Consistency of the estimator}
We aim to numerically assess the consistency of the estimator \(\hat{\theta}=(\lambda^{N,m},\nu^{N,m})\). To achieve this, we utilize the estimated values of \(\hat{\theta}=(\lambda^{N,m},\nu^{N,m})\) obtained from simulated datasets for various values of \(N\) and \(m\) to calculate the following empirical moments:
\begin{itemize}
	\item The empirical mean
	\begin{align*}
	\bar{\hat{\theta}}=(\overline{\lambda^{N,m}},\overline{\nu^{N,m}})=\frac{1}{100}\bigg(\sum_{i=1}^{100}\lambda^{N,m}_i,\sum_{i=1}^{100}\nu^{N,m}_i\bigg)\approx\mathbb{E}(\lambda^{N,m},\nu^{N,m}).
	\end{align*}	
	\item The empirical standard deviation
	\begin{align*}
		s_{\hat{\theta}}=(s_{\lambda^{N,m}},s_{\nu^{N,m}})=\sqrt{\frac{1}{99}\bigg(\sum_{i=1}^{100}(\lambda^{N,m}_i-\overline{\lambda^{N,m}})^2,\sum_{i=1}^{100}(\nu^{N,m}_i-\overline{\nu^{N,m}})^2\bigg)}\approx\sqrt{\mathbb{V}(\lambda^{N,m},\nu^{N,m})}.
	\end{align*}
	\item The empirical mean square error
	\begin{align*}
		MSE=\frac{1}{100}\bigg(\sum_{i=1}^{100}(\lambda^{N,m}_i-\lambda)^2,\sum_{i=1}^{100}(\nu^{N,m}_i-\nu)^2\bigg)\approx\mathbb{E}\big((\lambda^{N,m}-\lambda)^2,(\nu^{N,m}-\nu)^2\big).
	\end{align*}
	\item The empirical mean error
	 	\begin{align*}
	 	{\rm Mean}-{\rm Error}=\frac{1}{100}\bigg(\sum_{i=1}^{100}(\lambda^{N,m}_i-\lambda),\sum_{i=1}^{100}(\nu^{N,m}_i-\nu)\bigg)\approx\mathbb{E}\big((\lambda^{N,m}-\lambda),(\nu^{N,m}-\nu)\big).
	 \end{align*}
\end{itemize}
The results are presented in Table \ref{moment-est}. As observed, as \(N\) and \(m\) increase, the empirical mean of the estimator converges to the true parameter values, while the corresponding standard deviations decrease. Furthermore, both the mean squared error and the absolute mean error diminish as \(N\) and \(m\) grow larger. These findings validate the consistency of the estimator, as stated in Theorem \ref{cons-theo}.
 \begin{table}[h]
 	\begin{center}
 		\begin{tabular}{|c|c|c|c|c|}
 			\hline
 			Parameters moments & $\bar{\hat{\theta}}=(\overline{\lambda^{N,m}},\overline{\nu^{N,m}})$ & $s_{\hat{\theta}}=(s_{\lambda^{N,m}},s_{\nu^{N,m}})$  & MSE & Mean Error \\
 			\hline
 			$N=100$, $m=1000$ & $(-0.03, 0.33)$ & $(0.13, 0.14)$ &  $(0.30, 0.45)$  & $(-0.53, -0.66)$ \\
 			\hline
 		$N=500$, $m=10000$	& $(0.18,0.61)$  & $(0.13, 0.15)$ & $(0.12, 0.17)$  & $(-0.31, -0.38)$  \\
 			\hline
 		  $N=1000$, $m=10000$  & $0.28,0.74)$ & $(0.10, 0.11)$  & $(0.05, 0.07)$  & $(-0.21, -0.25)$ \\
 			\hline
 		$N=2000$, $m=20000$	 & $(0.46, 0.95)$ & $(0.08, 0.08)$ & $(0.008, 0.009)$    &  $(-0.03, -0.04)$\\
 			\hline
 			$N=3000$, $m=30000$	 & $(0.48, 0.97)$ & $(0.07, 0.07)$ & $(0.005, 0.006)$    &  $(-0.01, -0.02)$\\
 			\hline
 		\end{tabular}
 	\end{center}	
 	\caption{Empirical moments of the approximate LSE $\theta^{N,m}=(\lambda^{N,m},\nu^{N,m})$}
 	\label{moment-est}
 \end{table}

\subsubsection{Asymptotic normality}
In this section, we assess the asymptotic normality of the approximate LSE $\hat{\theta} = (\lambda^{N,m}, \nu^{N,m})$ as established in Theorem \ref{thm9}. Specifically, we numerically verify that the normalized error term $\sqrt{N}\big((\lambda^{N,m}-\lambda), (\nu^{N,m}-\nu)\big)$ converges to a Gaussian distribution as $N$ and $m$ approach infinity.

We utilize the simulated 100 samples from the power of two choices model with the true parameter $\theta = (\lambda, \nu) = (0.5, 1)$, varying the network size $N$ and the number of observations $m$. For each combination of $(N,m)$, we compute the first four empirical moments of the normalized error terms $\sqrt{N}\big((\lambda^{N,m}-\lambda), (\nu^{N,m}-\nu)\big)$: the mean, variance, skewness, and kurtosis. The results are summarized in Table \ref{mom-err}. Notably, we observe that the skewness and kurtosis values tend to approximate those of a normal distribution (0 for skewness and 3 for kurtosis), even for relatively small values of network size $N$ and number of observations $m$.
\begin{table}[h]
	\begin{center}
		\begin{tabular}{|c|c|c|c|c|}
			\hline
			Normalized error moments & 	Mean & Variance & Skewness & 	Kurtosis\\
			\hline
		N=100,m=1000 & $(-5.37, -6.60)$ &  $(1.75, 2.07)$  & $(-0.16, -0.41)$  & $(3.44, 4.31)$ \\
			\hline
		N=500,m=10000 & $(-7.15, -8.65)$ & $(8.99, 11.32)$ & $(0.51, 0.60)$ & $(3.90,4.21)$\\
			\hline
		N=1000,m=10000	 & $(-6.89, -8.08)$  & $(11.08, 14.38)$  & $(0.01, 0.02)$ & $(2.72, 2.51)$ \\
			\hline
		N=2000, m=20000  & $(-1.47, -2.01)$  &  $(13.96, 15.59)$  & $(0.07, 0.20)$ & $(2.80,2.68)$\\
			\hline
		N=3000, m=30000  & $(-0.86, -1.48)$  &  $(17.00, 17.15)$  & $(0.10, 0.26)$ & $(3.40, 3.07)$\\
			\hline	
		\end{tabular}
	\end{center}
	\caption{Empirical moments of the normalized errors $\sqrt{N}\big((\lambda^{N,m}-\lambda),(\nu^{N,m}-\nu)\big)$ }	
	\label{mom-err}
\end{table}

To further substantiate our findings, we test the normality of the normalized error terms using a Kolmogorov-Smirnov test. We conduct this test on the $100$ simulated datasets across the different values of network size $N$ and number of observations $m$. The resulting $p$-values are presented in Table \ref{kolmo-test}. As shown, the $p$-values are sufficiently large, suggesting that the null hypothesis asserting that the error terms $\sqrt{N}\big((\lambda^{N,m}-\lambda), (\nu^{N,m}-\nu)\big)$ follow a normal distribution is not rejected, even for smaller values of $N$ and $m$. This result indicates that the error terms tend to the normal distribution quickly. Nevertheless, it is also noted that the $p$-values are high for all $(N,m)$ combinations, and they do not necessarily increase with larger network sizes \(N\) and higher numbers of observations \(m\). This effect may be attributed to the fact that the true mean and variance of the error terms are unknown, requiring the use of empirical values, which may account for the observed outcomes.
\begin{table}[h]
	\begin{center}
	\begin{tabular}{|c|c|c|}
		\hline
	Network and data sizes $(N,m)$ & P-value for $ \sqrt{N}(\lambda^{N,m}-\lambda)$ & P-value for $ \sqrt{N}(\nu^{N,m}-\nu)$ \\
		\hline	
	$(N,m)=(100,1000)$  & $0.71$ & $0.48$ \\
				\hline
		$(N,m)=(500,10000)$  & $0.55$  & $0.59$ \\
		\hline
		$(N,m)=(1000,10000)$  & $0.99$ & $0.33$\\
		\hline
		 $(N,m)=(2000,20000)$  & $0.97$ & $0.94$ \\
		   \hline
		 $(N,m)=(3000,30000)$  & $0.65$ & $0.66$ \\
		   \hline
	\end{tabular}
	\end{center}
	\caption{Kolmogorov-Smirnov tests for the normalized parameter estimator errors $\sqrt{N}\big((\lambda^{N,m}-\lambda),(\nu^{N,m}-\nu)\big)$}
	\label{kolmo-test}
\end{table}

\begin{figure}[H]
	\begin{center}
	\includegraphics[width=0.4\textwidth]{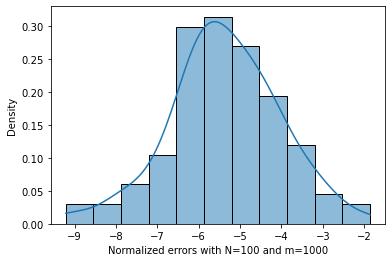}
	\includegraphics[width=0.4\textwidth]{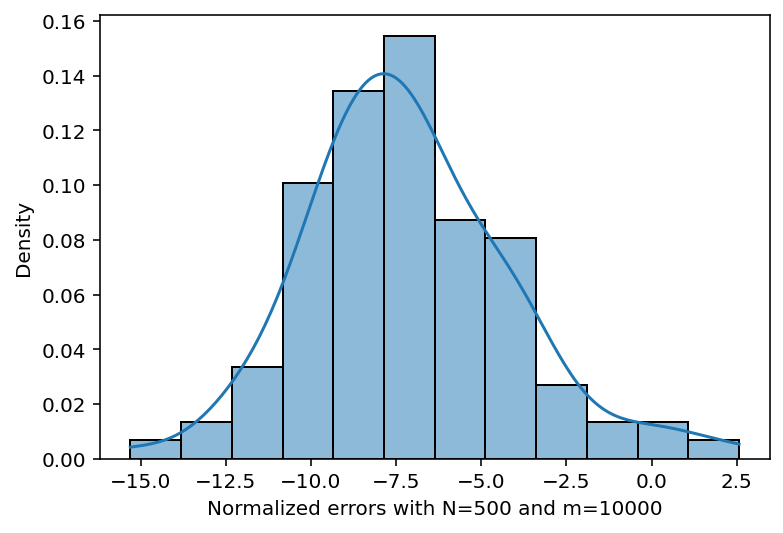}
	\includegraphics[width=0.4\textwidth]{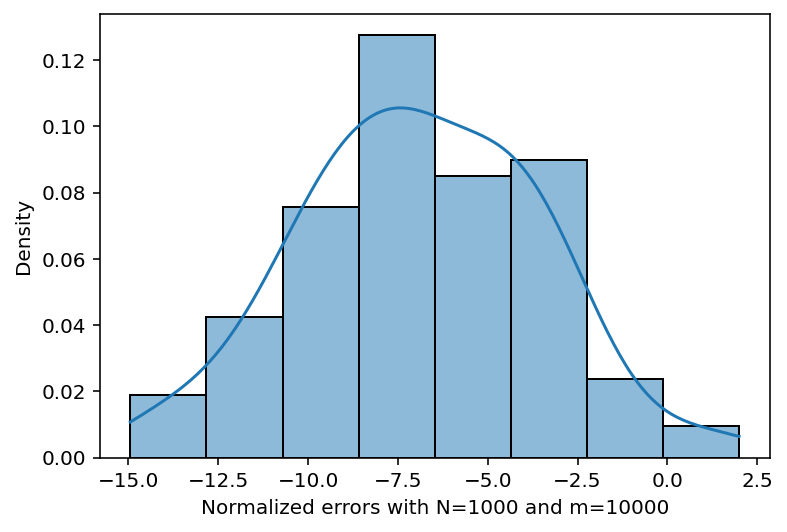}	
	\includegraphics[width=0.4\textwidth]{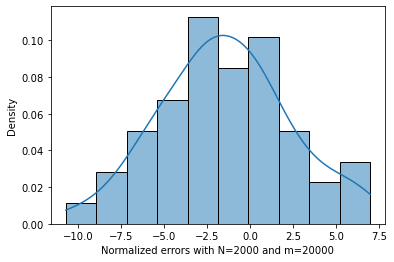}
	\includegraphics[width=0.4\textwidth]{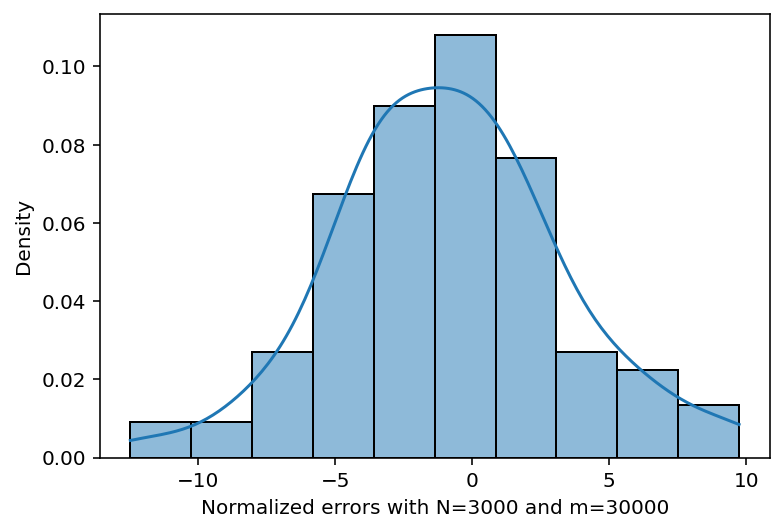}
	\end{center}
	\caption{Histograms of the normalized errors $\sqrt{N}(\lambda^{N,m}-\lambda)$ with the associated kernel density estimate plots for different values of the network size $N$ and the number of observations $m$ with the true parameters $\theta=(\lambda,\nu)=(0.5,1)$}
	\label{hist_lamb}
\end{figure}

\begin{figure}[H]
	\begin{center}
	\includegraphics[width=0.4\textwidth]{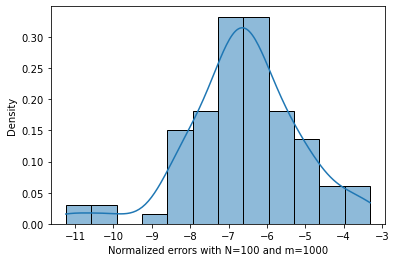}
	\includegraphics[width=0.4\textwidth]{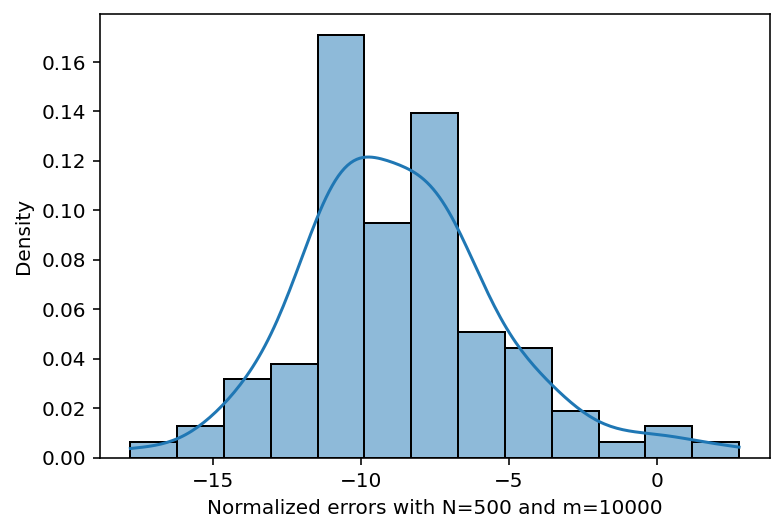}
		\includegraphics[width=0.4\textwidth]{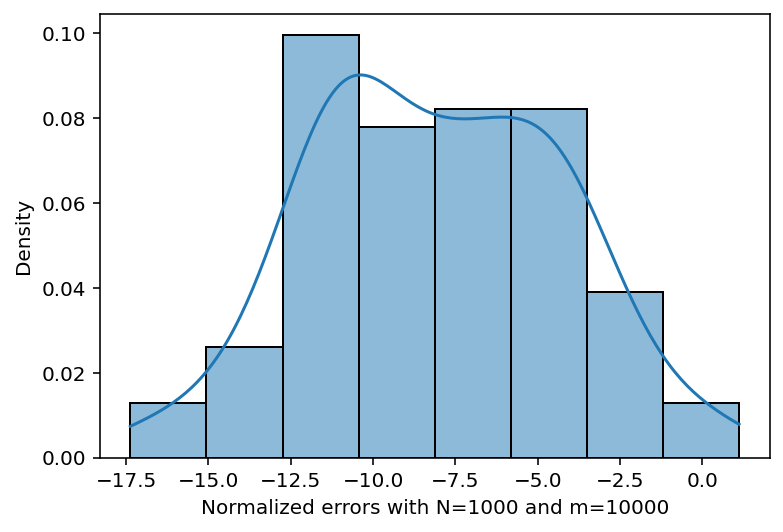}	
		\includegraphics[width=0.4\textwidth]{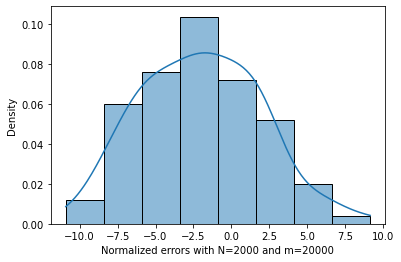}
		\includegraphics[width=0.4\textwidth]{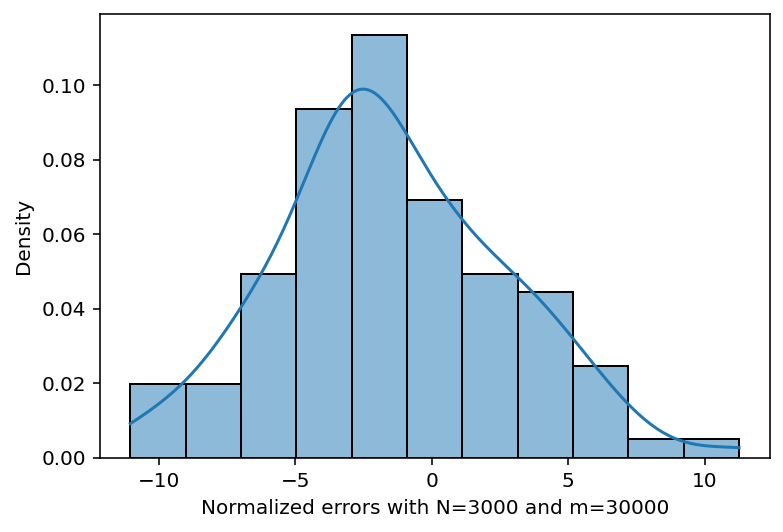}
	\end{center}
	\caption{Histograms of the normalized errors $\sqrt{N}(\nu^{N,m}-\nu)$ with the associated kernel density estimate plots for different values of the network size $N$ and the number of observations $m$, corresponding to true parameters $\theta=(\lambda,\nu)=(0.5,1)$}
	\label{hist_nu}
\end{figure}

Finally, we plot the histograms of the normalized error terms \(\sqrt{N}\big((\lambda^{N,m}-\lambda)\big)\) and \(\sqrt{N}\big((\nu^{N,m}-\nu)\big)\), along with a kernel density estimator. The results are shown in Figure \ref{hist_lamb} and Figure \ref{hist_nu}. Once again, we observe that the assertion of normality for the error terms aligns well with the empirical data.

\section{Conclusions and perspectives}
\label{conc}
In this paper, we considered the parameter estimation problem of the supermarket model. Based on an aggregate dataset, we constructed an approximate LSE by exploiting the law of large numbers together with the central limit theorem established for the model in the literature. Moreover, we proved the consistency together with the asymptotic normality of the estimator as both the size of the network and the number of observations go to infinity. Finally, we presented a numerical study where we tested our estimator against synthetic data obtained by simulating the power-of-two model highlighting our theoretical results.

The current work is the first statistical scheme for mean-field queuing systems and opens a new perspective. One naturally aims to investigate the statistical inference problem for other models. For instance, one can investigate the approximate LSE approach to the model proposed in \cite{Budhi+Fried2019} for load balancing mechanisms in cloud storage systems which is a generalization of the supermarket model. The established law of large numbers together with the central limit theorem make the approximate LSE approach used in the current work conceivable, provided that one can overcome the technical difficulty arising from the more complicated mean-field limiting equation. Another variation of the supermarket model for which one can study the parameter estimation problem is the one introduced in \cite{Budh+Mukh+Wu2019} in which the servers can communicate with their neighbors and where the neighborhood relationships are described in terms of a suitable graph. Again, the limit as the number of servers goes to infinity was identified, which can be exploited to build a statistical scheme, however, no central limit theorem nor a stationary distribution was established, therefore, the asymptotic normality of the estimator cannot be obtained by the similar scheme used in the current paper. Another interesting open problem is the nonparametric estimation of the interaction kernel in general mean-field queuing systems studied in \cite{Daw+Zhao2005}. Indeed, one can consider the exploitation of the limiting mean-field equation to build an estimator. However, contrary to our current proposal where the unknown parameters enter linearly in the mean-field limiting equation, in the nonparametric estimation one needs to deal with an optimization problem in function space. A potential avenue is to exploit the stationary distribution to build an estimator in the stationary regime and then investigate a justification for the interchange of limits $N\rightarrow\infty$ and $t\rightarrow\infty$.

Finally, the problem of statistical inference for general mean-field models on discrete space remains open and, as mentioned in the introduction, very few references exist.

\appendix

\section{Proof of Propositoin~\ref{explicit-sol}}
\label{explicit-sol-proof}

 Let $\xi(t)=(\xi_j(t))^T$ be the solution to the following infinite-dimensional ODE:
$$
d\xi(t)=G\left(M(t)+\xi(t),\varrho(t)\right)dt,\quad \xi(0)=z_0\in \tilde{\ell}_2.
$$
Define ${\cal Z}(t):=M(t)+\xi(t)$. Then, ${\cal Z}(t)$ satisfies the SDE:
$$
d{\cal Z}(t)=G({\cal Z}(t),\varrho(t))dt+a(t)dW(t),\quad {\cal Z}(0)=z_0.
$$
By \cite[Proposition 2]{Budhi+Fried2019}, we know that ${\cal Z}(t)\in\tilde{\ell}_2$ for all $t\in[0,T]$ almost surely. By the estimate \cite[Page 69, line -1]{Budhi+Fried2019} and Fatou's lemma (cf.  \cite[Page 78, Proof of Theorem 2]{Budhi+Fried2019}), we get
$$
E\Bigg[\sup_{t\in [0, T]}\sum_{m=0}^{\infty}(m+1)^2M_m^2(t)\Bigg]<\infty.
$$
Then, $M(t), \xi(t)\in \ell_1$ for all $t\in[0,T]$ almost surely. Moreover, for any $j\in\mathbb{Z}_+$, we have $\xi_j(0)=(z_0)_j$, and

{\small\begin{eqnarray}\label{2}
\frac{d\xi_j(t)}{dt}&=&G_j({\cal Z}(t),\varrho(t))\nonumber\\
&=&2\lambda\sum_{m=j}^{\infty}\{[M_{j-1}(t)+\xi_{j-1}(t)]\varrho_m(t)+[M_{m}(t)+\xi_{m}(t)]\varrho_{j-1}(t)\nonumber\\
&&\ \ \ \ \ \ \ \ \ \ -[M_{j}(t)+\xi_{j}(t)]\varrho_{m+1}(t)-[M_{m+1}(t)+\xi_{m+1}(t)]\varrho_{j}(t)\}\nonumber\\
&&+\nu\{[M_{j+1}(t)+\xi_{j+1}(t)]-[M_{j}(t)+\xi_{j}(t)] \}\nonumber\\
&=&2\lambda [M_{j-1}(t)-M_j(t)]\sum_{m=j}^{\infty}\varrho_m(t)+2\lambda [\xi_{j-1}(t)-\xi_j(t)]\sum_{m=j}^{\infty}\varrho_m(t)\nonumber\\
&&+2\lambda[\varrho_{j-1}(t)-\varrho_{j}(t)]\sum_{m=j}^{\infty}[M_{m}(t)+\xi_{m}(t)]+2\lambda M_{j}(t)\varrho_{j}(t)+2\lambda \xi_{j}(t)\varrho_{j}(t)+2\lambda\varrho_{j}(t)[M_{j}(t)+\xi_{j}(t)]\nonumber\\
&&+\nu\{[M_{j+1}(t)-M_{j}(t)]+[\xi_{j+1}(t)-\xi_{j}(t)] \}\nonumber\\
&=&2\lambda [M_{j-1}(t)-M_j(t)]\sum_{m=j}^{\infty}\varrho_m(t)+2\lambda[\varrho_{j-1}(t)-\varrho_{j}(t)]\sum_{m=j}^{\infty}M_{m}(t)+\nu[M_{j+1}(t)-M_{j}(t)]\nonumber\\
&&+4\lambda M_{j}(t)\varrho_{j}(t)+4\lambda \xi_{j}(t)\varrho_{j}(t)\nonumber\\
&&+2\lambda[\varrho_{j-1}(t)-\varrho_{j}(t)]\sum_{m=j}^{\infty}\xi_{m}(t)+2\lambda [\xi_{j-1}(t)-\xi_j(t)]\sum_{m=j}^{\infty}\varrho_m(t)+\nu[\xi_{j+1}(t)-\xi_{j}(t)] \nonumber\\
&=&2\lambda \left(\sum_{m=j}^{\infty}\varrho_m(t)\right)[M_{j-1}(t)-M_j(t)]+2\lambda[\varrho_{j-1}(t)-\varrho_{j}(t)]\sum_{m=j}^{\infty}M_{m}(t)+4\lambda \varrho_{j}(t)M_{j}(t)+\nu[M_{j+1}(t)-M_{j}(t)]\nonumber\\
&&+2\lambda[\varrho_{j-1}(t)-\varrho_{j}(t)]\sum_{m=j+2}^{\infty}\xi_{m}(t)+\{2\lambda[\varrho_{j-1}(t)-\varrho_{j}(t)]+\nu\}\xi_{j+1}(t)\nonumber\\
&&+\left\{2\lambda\left[\varrho_{j-1}(t)+\varrho_{j}(t)-\sum_{m=j}^{\infty}\varrho_m(t)\right]-\nu\right\}\xi_{j}(t)+2\lambda\left(\sum_{m=j}^{\infty}\varrho_m(t)\right)\xi_{j-1}(t).
\end{eqnarray}}

\noindent Since $\sum_{m=0}^{\infty}{\cal Z}_{m}(t)=0$, we get
$$
\sum_{m=j+2}^{\infty}\xi_{m}(t)=-\sum_{m=0}^{\infty}M_m(t)-\sum_{m=0}^{j+1}\xi_{m}(t).
$$
Similarly, by $\sum_{m=0}^{\infty}\varrho_m(t)=1$, we get
$$
\sum_{m=j}^{\infty}\varrho_{m}(t)=1-\sum_{m=0}^{j-1}\varrho_m(t).
$$
Then, we can rewrite (\ref{2}) as follows:
{\small\begin{eqnarray}\label{3}
\frac{d\xi_j(t)}{dt}&=&2\lambda \left[1-\sum_{m=0}^{j-1}\varrho_m(t)\right][M_{j-1}(t)-M_j(t)]-2\lambda[\varrho_{j-1}(t)-\varrho_{j}(t)]\sum_{m=0}^{j-1}M_{m}(t)+4\lambda \varrho_{j}(t)M_{j}(t)+\nu[M_{j+1}(t)-M_{j}(t)]\nonumber\\
&&-2\lambda[\varrho_{j-1}(t)-\varrho_{j}(t)]\sum_{m=0}^{j-2}\xi_{m}(t)+2\lambda\left[1-\sum_{m=0}^{j-2}\varrho_m(t)-2\varrho_{j-1}(t)+\varrho_{j}(t)\right]\xi_{j-1}(t)\nonumber\\
&&+\left\{2\lambda\left[-1+\sum_{m=0}^{j-1}\varrho_m(t)+2\varrho_{j}(t)\right]-\nu\right\}\xi_{j}(t)+\nu\xi_{j+1}(t).
\end{eqnarray}}

\noindent Note that (\ref{3}) can be regarded as an infinite-dimensional non-autonomous linear system of ODEs with random coefficients. Define $g(t)=(g_j(t))_{j=0}^{\infty}$ and $C(t)=(C_{jl}(t))_{j,l=0}^{\infty}$ by
\begin{eqnarray}\label{4}
g_j(t)&:=&2\lambda \left[1-\sum_{m=0}^{j-1}\varrho_m(t)\right][M_{j-1}(t)-M_j(t)]-2\lambda[\varrho_{j-1}(t)-\varrho_{j}(t)]\sum_{m=0}^{j-1}M_{m}(t)\nonumber\\
&&+4\lambda \varrho_{j}(t)M_{j}(t)+\nu[M_{j+1}(t)-M_{j}(t)]\nonumber\\
&=&-2\lambda[\varrho_{j-1}(t)-\varrho_{j}(t)]\sum_{m=0}^{j-2}M_{m}(t)+2\lambda \left[1-\sum_{m=0}^{j-2}\varrho_m(t)-2\varrho_{j-1}(t)+\varrho_{j}(t)\right]M_{j-1}(t)\nonumber\\
&&+\left\{2\lambda\left[-1+\sum_{m=0}^{j-1}\varrho_m(t)+2\varrho_{j}(t)\right]-\nu\right\}M_j(t)+\nu M_{j+1}(t),
\end{eqnarray}
and
\begin{eqnarray}\label{444}
C_{jl}(t):=\left\{ \begin{array}{ll}
         -2\lambda[\varrho_{j-1}(t)-\varrho_{j}(t)],\ \  & \mbox{if $0\le l\le j-2$},\\
2\lambda\left[1-\sum_{m=0}^{j-2}\varrho_m(t)-2\varrho_{j-1}(t)+\varrho_{j}(t)\right],\ \  & \mbox{if $l=j-1$},\\
        2\lambda\left[-1+\sum_{m=0}^{j-1}\varrho_m(t)+2\varrho_{j}(t)\right]-\nu,\ \  & \mbox{if $l=j$},\\
 \nu,\ \  & \mbox{if $l=j+1$},\\
 0,\ \  & \mbox{if $l\ge j+2$}.\end{array} \right.
\end{eqnarray}
Then, (\ref{3}) becomes
$$
\frac{d\xi(t)}{dt}=C(t)\xi(t)+g(t),\quad \xi(0)=z_0.
$$

Denote ${\cal C}=6\lambda+\nu.$ Therefore, for any $j\in\mathbb{Z}_+$ and $t\in[0,T]$, by (\ref{4}), we get
\begin{eqnarray*}
|g_j(t)|\le{\cal C}\Bigg\{|M_{j-1}(t)|+|M_j(t)|+|M_{j+1}(t)|+[\varrho_{j-1}(t)+\varrho_{j}(t)]\Bigg[\frac{\pi^2}{6}\sup_{t\in [0,T]}\sum_{m=0}^{\infty}(m+1)^2M^2_{m}(t)\Bigg]^{\frac{1}{2}}\Bigg\}.
\end{eqnarray*}
Then,
\begin{eqnarray*}
\sup_{t\in [0,T]}\|g(t)\|_1\le5{\cal C}\Bigg[\frac{\pi^2}{6}\sup_{t\in [0,T]}\sum_{m=0}^{\infty}(m+1)^2M^2_{m}(t)\Bigg]^{\frac{1}{2}}.
\end{eqnarray*}
For $x\in \ell_1$, $j\in\mathbb{Z}_+$ and $t\in[0,T]$, by (\ref{444}), we get
\begin{eqnarray*}
|(C(t)x)_j|\le{\cal C}\Bigg\{[\varrho_{j-1}(t)+\varrho_{j}(t)]\sum_{l=0}^{j-2}|x_l|+|x_{j-1}|+|x_j|+|x_{j+1}|\Bigg\}.
\end{eqnarray*}
Then,
\begin{eqnarray*}
\|C(t)x\|_1\le5{\cal C}\|x\|_1.
\end{eqnarray*}
Hence, by induction, we obtain that
$$
\|[C(t)]^nx\|_1\le(5{\cal C})^n\|x\|_1,\quad n\in \mathbb{N}.
$$
Thus, we have the following explicit expressions:
\begin{eqnarray}\label{5}
\xi(t)=e^{\int_0^tC(s)ds}z_0+\int_0^te^{\int_s^tC(u)du}g(s)ds,\ \ \ \ {\cal Z}(t)=M(t)+\xi(t).
\end{eqnarray}
Since $M(t)$ is a Gaussian martingale, by (\ref{4}), we deduce that the distributions of $\xi(t)$ and ${\cal Z}(t)$ are both Gaussian. The proof is complete.

\section{Proof of Lemma~\ref{lem:3.1}}
\label{lem3-1-proof}

 By H\"older's inequality, we find that $a_{11}a_{22}\ge(a_{12})^2$ and the equality sign holds if and only if $a_{11}=0$, or $a_{22}=0$, or
$$
a_{11}, a_{22}>0\ \ \ \ {\rm and}\ \ \ \
$$
for all  $j'\in\mathbb{Z}_+$ (see, e.g. \cite{Steele2004})
\begin{equation}
	\begin{split}
		&  \frac{\int_0^T U_{j'}(\varrho(s,\theta^*))ds}{\bigg(\sum\limits_{j=0}^{\infty}\left[\int_0^T U_j(\varrho(s,\theta^*))ds\right]^2\bigg)^{\frac{1}{2}}}=\frac{\int_0^T V_{j'}(\varrho(s,\theta^*))ds}{\bigg(\sum\limits_{j=0}^{\infty}\left[\int_0^T V_j(\varrho(s,\theta^*))ds\right]^2\bigg)^{\frac{1}{2}}}.
			\end{split}
\label{Hold-ineq}
\end{equation}
Note that for $s\in[0,T]$,
\begin{eqnarray*}
	\begin{split}
		U_0(\varrho(s,\theta^*))= -2 \varrho_0(s,\theta^*) \sum_{i=1}^{\infty}\varrho_i(s,\theta^*) -(\varrho_0(s,\theta^*))^2 ,
	\end{split}
\end{eqnarray*}
and
\begin{eqnarray*}
	\begin{split}
		U_1(\varrho(s,\theta^*))= 2 \varrho_0(s,\theta^*) \sum_{i=1}^{\infty}\varrho_i(s,\theta^*)  -2 \varrho_1(s,\theta^*) \sum_{i=2}^{\infty}\varrho_i(s,\theta^*)+(\varrho_0(s,\theta^*))^2-(\varrho_1(s,\theta^*))^2 ,
	\end{split}
\end{eqnarray*}
Then, by the fact that $\sum\limits_{i=0}^{\infty}\varrho_i(s,\theta^*)=1$, we get
\begin{eqnarray*}
	\begin{split}
		U_0(\varrho(s,\theta^*))&= -2 \varrho_0(s,\theta^*)  (1-\varrho_0(s,\theta^*))      -(\varrho_0(s,\theta^*))^2 =   \varrho_0(s,\theta^*)(  \varrho_0(s,\theta^*)-2),
	\end{split}
\end{eqnarray*}
and
\begin{eqnarray*}
	\begin{split}
		U_1(\varrho(s,\theta^*))&= 2 \varrho_0(s,\theta^*)  (1-\varrho_0(s,\theta^*))     -2 \varrho_1(s,\theta^*)  (1-\varrho_0(s,\theta^*)-\varrho_1(s,\theta^*)) +(\varrho_0(s,\theta^*))^2-(\varrho_1(s,\theta^*))^2\\
		&=   (\varrho_0(s,\theta^*)-\varrho_1(s,\theta^*))( 2- \varrho_0(s,\theta^*))+(\varrho_1(s,\theta^*))(\varrho_0(s,\theta^*)+\varrho_1(s,\theta^*)).
	\end{split}
\end{eqnarray*}
We have
$$
V_0(\varrho(s,\theta^*))=\varrho_1(s,\theta^*)-\varrho_0(s,\theta^*),
$$
and
$$
V_1(\varrho(s,\theta^*))=\varrho_2(s,\theta^*)-\varrho_1(s,\theta^*).
$$
\vskip 0.2cm
\noindent (a) Suppose that $(\ref{cond-ineq-posed1})$ holds. Then, $\int_0^TU_0(\varrho(s,\theta^*))ds< 0$ and $\int_0^TV_0(\varrho(s,\theta^*))ds> 0$. Thus, $(\ref{Hold-ineq})$ cannot hold  and hence  $(\ref{equal-Hold})$ holds.
\vskip 0.2cm
\noindent (b) Suppose that $(\ref{cond-ineq-posed2})$ holds. Then, $\int_0^TU_1(\varrho(s,\theta^*))ds> 0$ and $\int_0^TV_1(\varrho(s,\theta^*))ds< 0$. Thus, $(\ref{Hold-ineq})$ cannot hold  and hence  $(\ref{equal-Hold})$ holds.

\section{Proof of Theorem \ref{cons-theo}}
\label{theo-3-1-proof}
By Theorem \ref{LLN-Budhi}, \cite[(5.7), page 117  and Proposition 5.3, page 119]{Eth+Kur86} and Lemma \ref{lem1}, we obtain that
\begin{equation}
	\sup_{t\in[0, T]}\sum_{j=0}^{\infty}|\varrho^N_j(t,\theta^*)-\varrho_j(t,\theta^*)|\rightarrow 0
	\label{Mar8a}
\end{equation} 
 in probability as $N\rightarrow\infty$. Moreover, 
\begin{equation*}
	\begin{split}
	|a^{N,m}_{11}-a_{11}|&\leq\sum\limits_{j=0}^{\infty}\left|\left[\frac{T}{m}\sum\limits_{k=1}^mU_j(\varrho^N(t_k,\theta^*))\right]^2-\left[\int_0^T U_j(\varrho(s))ds\right]^2\right|\\
	&\leq {\sup_{j\in\mathbb{Z}_+}}\bigg\{\left|\frac{T}{m}\sum\limits_{k=1}^mU_j(\varrho^N(t_k,\theta^*))+\int_0^T U_j(\varrho(s))ds\right|\bigg\}\\
	&\qquad\times\sum\limits_{j=0}^{\infty}\left|\frac{T}{m}\sum\limits_{k=1}^mU_j(\varrho^N(t_k,\theta^*))-\int_0^T U_j(\varrho(s))ds\right|\\
	&\leq T{\sup_{j\in\mathbb{Z}_+}\Big\{\sup_{ t\in[0, T]}\big|U_j(\varrho^N(t,\theta^*))\big|+\sup_{t\in[0,T] }\big|U_j(\varrho(t,\theta^*))\big|\Big\}}\\
	&\qquad\times\sum\limits_{j=0}^{\infty}\bigg\{\left|\frac{T}{m}\sum\limits_{k=1}^mU_j(\varrho^N(t_k,\theta^*))-\frac{T}{m}\sum\limits_{k=1}^mU_j(\varrho(t_k,\theta^*))\right|\\
	&\qquad\quad\qquad+\left|\frac{T}{m}\sum\limits_{k=1}^mU_j(\varrho(t_k,\theta^*))-\int_0^T U_j(\varrho(s))ds\right|\bigg\}.
\end{split}
\end{equation*}
By (\ref{U-func}), we get $|U_j(x)|\leq 6\|x\|_1$ for all $x\in\ell_1$ and $j\in\mathbb{Z}_+$. Then
\begin{align*}
	|a^{N,m}_{11}-a_{11}|&\leq 12T\bigg\{T\sup_{t\in[0,T]}\sum\limits_{j=0}^{\infty}\left|U_j(\varrho^N(t,\theta^*))-U_j(\varrho(t,\theta^*))\right|\\
	&\qquad\quad+\sum\limits_{j=0}^{\infty}\left|\frac{T}{m}\sum\limits_{k=1}^mU_j(\varrho(t_k,\theta^*))-\int_0^T U_j(\varrho(s))ds\right|\bigg\}.
\end{align*}
Therefore, by (\ref{Mar8a}), we get that the right hand side of the last inequality goes to $ 0\ \ {\rm in\ probability}\ {\rm as}\ N,m\rightarrow\infty$. 
Similarly, we can show that
\begin{equation*}
	|a^{N,m}_{12}-a_{12}|,\,|a^{N,m}_{22}-a_{22}|,\,|b^{N,m}_{1}-b_{1}|,\,|b^{N,m}_{2}-b_{2}|\rightarrow 0\ \ {\rm in\ probability}\ {\rm as}\ N,m\rightarrow\infty.
\end{equation*}
Therefore, the proof is complete.

\section{Proof of Theorem~\ref{thm9}}
\label{theo-2-proof}
By $(\ref{lambda-nu-star})$ and $(\ref{Nov1b})$, we get
	\begin{equation*}
		\begin{split}
\sqrt{N}\bigg( \begin{pmatrix}
		\lambda^{N,m}\\
		\nu^{N,m}
	\end{pmatrix}- \begin{pmatrix}
		\lambda^*\\
		\nu^*
	\end{pmatrix}\bigg)=\sqrt{N} \begin{pmatrix}
		\frac{a^{N,m}_{22}b^{N,m}_1-a^{N,m}_{12}b^{N,m}_2}{a^{N,m}_{11}a^{N,m}_{22}-(a^{N,m}_{12})^2}-\frac{a_{22}b_1-a_{12}b_2}{a_{11}a_{22}-(a_{12})^2}\\
		\\
		\frac{-a^{N,m}_{21}b^{N,m}_1+a^{N,m}_{11}b^{N,m}_2 }{a^{N,m}_{11}a^{N,m}_{22}-(a^{N,m}_{12})^2}-\frac{-a_{21}b_1+a_{11}b_2 }{a_{11}a_{22}-(a_{12})^2}
	\end{pmatrix} .
\end{split}
\end{equation*}
Moreover, simple calculations lead to
\begin{equation*}
\begin{split}
\frac{a^{N,m}_{22}b^{N,m}_1-a^{N,m}_{12}b^{N,m}_2}{a^{N,m}_{11}a^{N,m}_{22}-(a^{N,m}_{12})^2}-\frac{a_{22}b_1-a_{12}b_2}{a_{11}a_{22}-(a_{12})^2}&=\frac{(a^{N,m}_{22}b^{N,m}_1-a^{N,m}_{12}b^{N,m}_2)-(a_{22}b_1-a_{12}b_2)}{a^{N,m}_{11}a^{N,m}_{22}-(a^{N,m}_{12})^2}\\
		&\quad+\frac{{(a_{22}b_1-a_{12}b_2)[(a_{11}a_{22}-(a_{12})^2)- (a^{N,m}_{11}a^{N,m}_{22}-(a^{N,m}_{12})^2) ]}}{(a^{N,m}_{11}a^{N,m}_{22}-(a^{N,m}_{12})^2)(a_{11}a_{22}-(a_{12})^2)},
\end{split}
\end{equation*}
and
\begin{equation*}
\begin{split}
\frac{-a^{N,m}_{21}b^{N,m}_1+a^{N,m}_{11}b^{N,m}_2 }{a^{N,m}_{11}a^{N,m}_{22}-(a^{N,m}_{12})^2}-\frac{-a_{21}b_1+a_{11}b_2 }{a_{11}a_{22}-(a_{12})^2}&=	\frac{(-a^{N,m}_{21}b^{N,m}_1+a^{N,m}_{11}b^{N,m}_2)-(-a_{21}b_1+a_{11}b_2)}{a^{N,m}_{11}a^{N,m}_{22}-(a^{N,m}_{12})^2}\\
&\quad+\frac{{(-a_{21}b_1+a_{11}b_2)[(a_{11}a_{22}-(a_{12})^2)- (a^{N,m}_{11}a^{N,m}_{22}-(a^{N,m}_{12})^2) ]}}{(a^{N,m}_{11}a^{N,m}_{22}-(a^{N,m}_{12})^2)(a_{11}a_{22}-(a_{12})^2)}.
\end{split}
\end{equation*}
To simplify notation, let us define
\begin{align*}
	\mathcal{I}^{N,m}:= \sqrt{N} \bigg[  (a^{N,m}_{22}b^{N,m}_1-a^{N,m}_{12}b^{N,m}_2) - (a_{22}b_1-a_{12}b_2) \bigg],
\end{align*}
\begin{align*}
	\mathcal{J}^{N,m}:=\sqrt{N} \bigg[   ( -a^{N,m}_{21}b^{N,m}_1+a^{N,m}_{11}b^{N,m}_2)-(-a_{21}b_1+a_{11}b_2)\bigg],
\end{align*}
\begin{align*}
	\mathcal{H}^{N,m}:=a^{N,m}_{11}a^{N,m}_{22}-(a^{N,m}_{12})^2,
\end{align*}
and
\begin{align*}
	\mathcal{K}^{N,m}:=\sqrt{N}\bigg[(a_{11}a_{22}-(a_{12})^2)- (a^{N,m}_{11}a^{N,m}_{22}-(a^{N,m}_{12})^2) \bigg].
\end{align*}
We will analyze the convergence of $\mathcal{I}^{N,m},\mathcal{J}^{N,m}$ and $\mathcal{K}^{N,m}$ as $N,m,\frac{m}{\sqrt{N}}\rightarrow\infty$. By the Skorohod representation theorem (cf. \cite[Page 102]{Eth+Kur86}), we can and do assume without loss of generality that $\mathcal{Z}^N$ converges to $\mathcal{Z}$ in probability in $\mathcal{D}([0, T ] , \ell_2)$. To save space, we will prove the convergence of $\mathcal{I}^{N,m}$, the convergence of the other terms follows by similar arguments.
\vskip 0.1cm
First, we have
\begin{align*}
	\mathcal{I}^{N,m}=\mathcal{I}^{N,m}_1+\mathcal{I}^{N,m}_2
\end{align*}
with
\begin{align*}
	\mathcal{I}^{N,m}_1=\sqrt{N}\big( a^{N,m}_{22}b^{N,m}_1 -a_{22}b_1\big)\ \ \mbox{ and }\ \
	\mathcal{I}^{N,m}_2=\sqrt{N}\big(a_{12}b_2-a^{N,m}_{12}b^{N,m}_2  \big).
\end{align*}
Moreover,
\begin{align*}
	\mathcal{I}^{N,m}_1= \sqrt{N} a_{22}(b_1^{N,m}-b_1)+\sqrt{N} b_1^{N,m}(a^{N,m}_{22}-a_{22}),
\end{align*}
and
\begin{align*}
	\mathcal{I}^{N,m}_2=-\sqrt{N}\big(a_{12}(b_2^{N,m}-b_2)+b_2^{N,m} (a_{12}^{N,m}-a_{12}) \big).
\end{align*}
Then, by adding and subtracting terms, we get
\begin{align}\label{Oct21jk}
	&b_1^{N,m}-b_1\nonumber\\
&= \sum\limits_{j=0}^{\infty}\left[{\varrho}_j(T,\theta^*)-{\varrho}_j(0,\theta^*)\right]\Bigg\{\left[\frac{T}{m}\sum\limits_{k=1}^mU_j(\varrho^N(t_k,\theta^*))\right]-\left[\frac{T}{m}\sum\limits_{k=1}^mU_j(\varrho(t_k,\theta^*))\right]\Bigg\}\nonumber\\
&\quad+ \sum\limits_{j=0}^{\infty}\Bigg\{\left[{\varrho}^N_j(T,\theta^*)-{\varrho}^N_j(0,\theta^*)\right]-\left[{\varrho}_j(T,\theta^*)-{\varrho}_j(0,\theta^*)\right]\Bigg\}\Bigg\{\left[\frac{T}{m}\sum\limits_{k=1}^mU_j(\varrho^N(t_k,\theta^*))\right]-\left[\frac{T}{m}\sum\limits_{k=1}^mU_j(\varrho(t_k,\theta^*))\right]\Bigg\}\nonumber\\
&\quad+ \sum\limits_{j=0}^{\infty}\left[{\varrho}^N_j(T,\theta^*)-{\varrho}^N_j(0,\theta^*)\right]\Bigg\{\left[\frac{T}{m}\sum\limits_{k=1}^mU_j(\varrho(t_k,\theta^*))\right]-\int_0^T U_j(\varrho(s))ds\Bigg\}\nonumber\\
	&\quad + \sum\limits_{j=0}^{\infty}\Bigg\{\left[{\varrho}^N_j(T,\theta^*)-{\varrho}^N_j(0,\theta^*)\right]-\left[{\varrho}_j(T,\theta^*)-{\varrho}_j(0,\theta^*)\right]\Bigg\}\int_0^T U_j(\varrho(s))ds.
\end{align}
Below we consider the convergence of each summation term in (\ref{Oct21jk}).
\vskip 0.1cm
\noindent (a) By Theorem $\ref{rrr3}$, as $N\rightarrow\infty$, the term
\begin{align*}
\sqrt{N}\sum\limits_{j=0}^{\infty}\Bigg\{\left[{\varrho}^N_j(T,\theta^*)-{\varrho}^N_j(0,\theta^*)\right]-\left[{\varrho}_j(T,\theta^*)-{\varrho}_j(0,\theta^*)\right]\Bigg\}\int_0^T U_j(\varrho(s))ds
\end{align*}
converges in probability to
\begin{align*}
	\sum\limits_{j=0}^{\infty}\left[\mathcal{Z}_j(T,\theta^*)-\mathcal{Z}_j(0,\theta^*)\right]\int_0^T U_j(\varrho(s))ds.
\end{align*}

\vskip 0.1cm
\noindent (b)
By (\ref{ODE-LLN-Budhi}) and (\ref{U-func})--(\ref{Oct21b}), we get
\begin{eqnarray*}
|U_j(\varrho(t_k,\theta^*))-U_j(\varrho(s))|&\le&2 \varrho_{j-1}(t_k,\theta^*)\sum_{i=0}^{\infty}\big|\varrho_{i}(t_k,\theta^*)-\varrho_{i}(s,\theta^*)\big|+2\big|\varrho_{j-1}(t_k,\theta^*)-\varrho_{j-1}(s,\theta^*)\big|\\
&&+2 \varrho_{j}(t_k,\theta^*)\sum_{i=0}^{\infty}\big|\varrho_{i}(t_k,\theta^*)-\varrho_{i}(s,\theta^*)\big|+2\big|\varrho_{j}(t_k,\theta^*)-\varrho_{j}(s,\theta^*)\big|\\
&&+2\big|\varrho_{j-1}(t_k,\theta^*)-\varrho_{j-1}(s,\theta^*)\big|+2\big|\varrho_{j}(t_k,\theta^*)-\varrho_{j}(s,\theta^*)\big|\\
&\le&6\sum_{i=0}^{\infty}\big|\varrho_{i}(t_k,\theta^*)-\varrho_{i}(s,\theta^*)\big|\\
&\le&6\sum_{i=0}^{\infty}\int_s^{t_k}\big|F_{i}(\varrho(u,\theta^*))\big|du\\
&\le&\frac{6(6\lambda+2\nu)T}{m}.
\end{eqnarray*}
Then,
\begin{equation*}
	\begin{split}
&\left|\sqrt{N}\sum\limits_{j=0}^{\infty}\left[{\varrho}^N_j(T,\theta^*)-{\varrho}^N_j(0,\theta^*)\right]\Bigg\{\left[\frac{T}{m}\sum\limits_{k=1}^mU_j(\varrho(t_k,\theta^*))\right]-\int_0^T U_j(\varrho(s))ds\Bigg\}\right|\\
&\qquad=\left|\sqrt{N}\sum\limits_{j=0}^{\infty}\left[{\varrho}^N_j(T,\theta^*)-{\varrho}^N_j(0,\theta^*)\right]\Bigg\{\sum\limits_{k=1}^m\int_{t_{k-1}}^{t_k} [U_j(\varrho(t_k,\theta^*))-U_j(\varrho(s))]ds\Bigg\}\right|\\
&\qquad\leq\sqrt{N}\sum\limits_{j=0}^{\infty}\left[{\varrho}^N_j(T,\theta^*)+{\varrho}^N_j(0,\theta^*)\right]\Bigg\{\sum\limits_{k=1}^m\int_{t_{k-1}}^{t_k} \frac{6(6\lambda+2\nu)T}{m}ds\Bigg\}\\
&\qquad=\frac{24(3\lambda+\nu)T^2\sqrt{N}}{m}\\
&\qquad \rightarrow0\ \ \ \ {\rm as}\ \frac{m}{\sqrt{N}}\rightarrow\infty.
\end{split}
\end{equation*}

\vskip 0.1cm
\noindent (c) By the estimate \cite[Page 70, line 4]{Budhi+Fried2019}, we have
 that
\begin{equation}\label{Mar8b}
\sup_{N\in\mathbb{N}}E\left[\sup_{t\in[0,T]}\sum_{j=0}^{\infty}(j+1)^2({\cal Z}^N_j(t,\theta^*))^2\right]<\infty.
\end{equation}
Then, by Fatou's lemma, we get
\begin{equation}\label{Mar8c}
E\left[\sup_{t\in[0,T]}\sum_{j=0}^{\infty}(j+1)^2({\cal Z}_j(t,\theta^*))^2\right]<\infty.
\end{equation}
By (\ref{ODE-LLN-Budhi}) and (\ref{U-func})--(\ref{Oct21b}), we get
\begin{eqnarray*}
|U_j(\varrho^N(t_k,\theta^*))-U_j(\varrho(t_k,\theta^*))|\le 6\sum_{i=0}^{\infty}\big|\varrho^N_{i}(t_k,\theta^*)-\varrho_{i}(t_k,\theta^*)\big|,
\end{eqnarray*}
which together with (\ref{Mar8a}) and (\ref{Mar8b}) implies that
\begin{eqnarray*}
&&\Bigg|\sqrt{N}\sum\limits_{j=0}^{\infty}\Bigg\{\left[{\varrho}^N_j(T,\theta^*)-{\varrho}^N_j(0,\theta^*)\right]-\left[{\varrho}_j(T,\theta^*)-{\varrho}_j(0,\theta^*)\right]\Bigg\}\\
&&\quad\cdot\Bigg\{\left[\frac{T}{m}\sum\limits_{k=1}^mU_j(\varrho^N(t_k,\theta^*))\right]-\left[\frac{T}{m}\sum\limits_{k=1}^mU_j(\varrho(t_k,\theta^*))\right]\Bigg\}\Bigg|\\
&=&\Bigg|\sum\limits_{j=0}^{\infty}\left[{\cal Z}^N_j(T,\theta^*)-{\cal Z}^N_j(0,\theta^*)\right]\Bigg\{\left[\frac{T}{m}\sum\limits_{k=1}^mU_j(\varrho^N(t_k,\theta^*))\right]-\left[\frac{T}{m}\sum\limits_{k=1}^mU_j(\varrho(t_k,\theta^*))\right]\Bigg\}\Bigg|\\
&\le&6\sum\limits_{j=0}^{\infty}\left[|{\cal Z}^N_j(T,\theta^*)|+|{\cal Z}^N_j(0,\theta^*)|\right]\Bigg\{\frac{T}{m}\sum\limits_{k=1}^m\sum_{i=0}^{\infty}|\varrho^N_i(t_k,\theta^*)-\varrho_i(t_k,\theta^*)|\Bigg\}\\
&\le&6T\sum\limits_{j=0}^{\infty}\left[|{\cal Z}^N_j(T,\theta^*)|+|{\cal Z}^N_j(0,\theta^*)|\right]\cdot\sup_{t\in[0,T]}\sum_{i=0}^{\infty}|\varrho^N_i(t,\theta^*)-\varrho_i(t,\theta^*)|\\
&\rightarrow&0\ \ {\rm in\ probability}\ {\rm as}\ N\rightarrow\infty.
\end{eqnarray*}

\vskip 0.1cm
\noindent (d) For $i\in \mathbb{Z}_+$ and $k\in\{1,2,\dots,m\}$, define a non-negative measure $\tau^{i,k}$ on $\mathbb{Z}_+$ by
\begin{eqnarray*}
	\tau^{i,k}_l:=\begin{cases}
		\varrho_l(t_k,\theta^*),&l< i , \\
		\varrho^N_l(t_k,\theta^*),&l\ge i .
	\end{cases}
	\end{eqnarray*}
Note that
\begin{eqnarray*}
	\partial_lU_j(x)=\begin{cases}
		0,&l<j-1, \\
2\sum_{p=j-1}^{\infty}x_p,&l=j-1, \\
2x_{j-1}-2\sum_{p=j}^{\infty}x_p,&l=j, \\
		2x_{j-1}-2x_j,&l\ge j+1.
	\end{cases}
	\end{eqnarray*}
Then, for $s\in[t_{k-1},t_k]$, we have

\begin{eqnarray}\label{Ocy21asd}
&&\left|\sqrt{N}[U_j(\varrho^N(t_k,\theta^*))-U_j(\varrho(t_k,\theta^*))]-\sum_{i=0}^{\infty}\partial_iU_j(\varrho(s))\mathcal{Z}^N_i(t_k,\theta^*)\right|\nonumber\\
&=&\left|\sum_{i=0}^{\infty}\left\{\sqrt{N}[U_j(\tau^{i,k})-U_j(\tau^{i+1,k})]-\partial_iU_j(\varrho(s))\mathcal{Z}^N_i(t_k,\theta^*)\right\}\right|\nonumber\\
&\le&\sum_{i=0}^{\infty}\left|\sqrt{N}[U_j(\tau^{i,k})-U_j(\tau^{i+1,k})]-\partial_iU_j(\varrho(s))\mathcal{Z}^N_i(t_k,\theta^*)\right|\nonumber\\
&=&\sum_{i=0}^{\infty}\Bigg|\sqrt{N}\int_{\varrho_i(t_k,\theta^*)}^{\varrho^N_i(t_k,\theta^*)}\partial_iU_j(\varrho_0(t_k,\theta^*),\dots,\varrho_{i-1}(t_k,\theta^*),u,\varrho^N_{i+1}(t_k,\theta^*),\varrho^N_{i+2}(t_k,\theta^*),\dots)du\nonumber\\
&&\quad\quad-\sqrt{N}\int_{\varrho_i(t_k,\theta^*)}^{\varrho^N_i(t_k,\theta^*)}\partial_iU_j(\varrho(s))du\Bigg|\nonumber\\
&\le&2\left\{\sup_{t\in[0,T]}\sum_{i=0}^{\infty}|\mathcal{Z}^N_i(t,\theta^*)|\right\}\left\{\sup_{1\le k\le m}\sum_{i=0}^{\infty}\sup_{s,t\in[t_{k-1},t_k]}|\varrho_i(s,\theta^*)-\varrho_i(t,\theta^*)|+\sup_{t\in[0,T]}\sum_{i=0}^{\infty}|\varrho_i^N(t,\theta^*)-\varrho_i(t,\theta^*)|\right\}.\nonumber\\
&&
\end{eqnarray}

By  \cite[(5.7), page 117  and Proposition 5.3, page 119]{Eth+Kur86}, we have that
\begin{equation}\label{Mar8d}
\sup_{t\in [0,T]}\|\mathcal{Z}^N(t,\theta^*)-\mathcal{Z}(t,\theta^*)\|_2\rightarrow 0
\end{equation}
in probability as $N\rightarrow\infty$. Then, by (\ref{Mar8a}) and (\ref{Mar8b})--(\ref{Mar8d}), we get

\begin{eqnarray*}
&&\Bigg|\sqrt{N}\sum\limits_{j=0}^{\infty}\left[{\varrho}_j(T,\theta^*)-{\varrho}_j(0,\theta^*)\right]\Bigg\{\left[\frac{T}{m}\sum\limits_{k=1}^mU_j(\varrho^N(t_k,\theta^*))\right]-\left[\frac{T}{m}\sum\limits_{k=1}^mU_j(\varrho(t_k,\theta^*))\right]\Bigg\}\\
&&-\sum\limits_{j=0}^{\infty}\left[{\varrho}_j(T,\theta^*)-{\varrho}_j(0,\theta^*)\right]	\Bigg[\sum_{i=0}^{\infty}\int_0^T\partial_iU_j(\varrho(s))\mathcal{Z}_i(s,\theta^*)ds\Bigg]\Bigg|\\
&\le&\sum\limits_{j=0}^{\infty}\left[{\varrho}_j(T,\theta^*)+{\varrho}_j(0,\theta^*)\right]\Bigg\{\sum\limits_{k=1}^m\int_{t_{k-1}}^{t_k}\left|\sqrt{N}[U_j(\varrho^N(t_k,\theta^*))-U_j(\varrho(t_k,\theta^*))]-\sum_{i=0}^{\infty}\partial_iU_j(\varrho(s))\mathcal{Z}_i(s,\theta^*)\right|ds\Bigg\}\\
&\le&2\sup_{j\in\mathbb{Z}_+}\Bigg\{\sum\limits_{k=1}^m\int_{t_{k-1}}^{t_k}\left|\sqrt{N}[U_j(\varrho^N(t_k,\theta^*))-U_j(\varrho(t_k,\theta^*))]-\sum_{i=0}^{\infty}\partial_iU_j(\varrho(s))\mathcal{Z}_i(s,\theta^*)\right|ds\Bigg\}\\
&\le&2\sup_{j\in\mathbb{Z}_+}\Bigg\{\sum\limits_{k=1}^m\int_{t_{k-1}}^{t_k}\left|\sqrt{N}[U_j(\varrho^N(t_k,\theta^*))-U_j(\varrho(t_k,\theta^*))]-\sum_{i=0}^{\infty}\partial_iU_j(\varrho(s))\mathcal{Z}^N_i(t_k,\theta^*)\right|ds\Bigg\}\\
&&+2\sup_{j\in\mathbb{Z}_+}\Bigg\{\sum\limits_{k=1}^m\int_{t_{k-1}}^{t_k}\left|\sum_{i=0}^{\infty}\partial_iU_j(\varrho(s))\mathcal{Z}^N_i(t_k,\theta^*)-\sum_{i=0}^{\infty}\partial_iU_j(\varrho(s))\mathcal{Z}_i(t_k,\theta^*)\right|ds\Bigg\}\\
&&+2\sup_{j\in\mathbb{Z}_+}\Bigg\{\sum\limits_{k=1}^m\int_{t_{k-1}}^{t_k}\left|\sum_{i=0}^{\infty}\partial_iU_j(\varrho(s))\mathcal{Z}_i(t_k,\theta^*)-\sum_{i=0}^{\infty}\partial_iU_j(\varrho(s))\mathcal{Z}_i(s,\theta^*)\right|ds\Bigg\}\\
&\le&4T\Bigg\{\sup_{t\in[0,T]}\sum_{i=0}^{\infty}|\mathcal{Z}^N_i(t,\theta^*)|\Bigg\}\Bigg\{\sup_{1\le k\le m}\sum_{i=0}^{\infty}\sup_{s,t\in[t_{k-1},t_k]}|\varrho_i(s,\theta^*)-\varrho_i(t,\theta^*)|+\sup_{t\in[0,T]}\sum_{i=0}^{\infty}|\varrho_i^N(t,\theta^*)-\varrho_i(t,\theta^*)|\Bigg\}\\
&&+12T\sup_{t\in[0,T]}\sum_{i=0}^{\infty}|\mathcal{Z}^N_i(t,\theta^*)-\mathcal{Z}_i(t,\theta^*)|+12\sum\limits_{k=1}^m\int_{t_{k-1}}^{t_k}\sum_{i=0}^{\infty}|\mathcal{Z}_i(t_k,\theta^*)-\mathcal{Z}_i(s,\theta^*)|ds\\
&\rightarrow&0\ \ {\rm in\ probability}\ {\rm as}\ N,m\rightarrow\infty.
\end{eqnarray*}
Thus, by (a)--(d), we deduce that as $N,m,\frac{m}{\sqrt{N}}\rightarrow\infty$, $\sqrt{N}(b_1^{N,m}-b_1)$  converges in probability to
	\begin{align}\label{N-b1nm-b1-lim}
	&\sum\limits_{j=0}^{\infty}\left[{\varrho}_j(T,\theta^*)-{\varrho}_j(0,\theta^*)\right]	\Bigg[\sum_{i=0}^{\infty}\int_0^T\partial_iU_j(\varrho(s))\mathcal{Z}_i(s,\theta^*)ds\Bigg]\nonumber\\
	&+\sum\limits_{j=0}^{\infty}\left[\mathcal{Z}_j(T,\theta^*)-\mathcal{Z}_j(0,\theta^*)\right]\int_0^T U_j(\varrho(s))ds.
	\end{align}

By adding and subtracting terms, we  get
\begin{align*}
	a_{12}^{N,m}-a_{12}&= \sum\limits_{j=0}^{\infty}\Bigg\{\left[\frac{T}{m}\sum\limits_{k=1}^mU_j(\varrho^N(t_k,\theta^*))\right]-\left[\frac{T}{m}\sum\limits_{k=1}^mU_j(\varrho(t_k,\theta^*))\right]\Bigg\}\left[\frac{T}{m}\sum\limits_{k=1}^mV_j(\varrho^N(t_k,\theta^*))\right]\\
	&\quad + \sum\limits_{j=0}^{\infty}\Bigg\{\left[\frac{T}{m}\sum\limits_{k=1}^mU_j(\varrho(t_k,\theta^*))\right]-\int_0^T U_j(\varrho(s))ds\Bigg\}\left[\frac{T}{m}\sum\limits_{k=1}^mV_j(\varrho^N(t_k,\theta^*))\right]\\
	&\quad + \sum\limits_{j=0}^{\infty}\int_0^T U_j(\varrho(s))ds\Bigg\{\left[\frac{T}{m}\sum\limits_{k=1}^mV_j(\varrho^N(t_k,\theta^*))\right]-\left[\frac{T}{m}\sum\limits_{k=1}^mV_j(\varrho(t_k,\theta^*))\right]\Bigg\}\\
&\quad + \sum\limits_{j=0}^{\infty}\int_0^T U_j(\varrho(s))ds\Bigg\{\left[\frac{T}{m}\sum\limits_{k=1}^mV_j(\varrho(t_k,\theta^*))\right]-\int_0^T V_j(\varrho(s))ds\Bigg\}.
\end{align*}
Similar to the above argument, we can show  that as  $N,m,\frac{m}{\sqrt{N}}\rightarrow\infty$, $\sqrt{N}(a_{12}^{N,m}-a_{12})$ converges in probability to
\begin{equation}
	\begin{split}
	&\sum\limits_{j=0}^{\infty}\int_0^TV_j(\varrho(s))ds \left[\sum_{i=0}^{\infty}\int_0^T\partial_iU_j(\varrho(s))\mathcal{Z}_i(s,\theta^*)ds\right]\\
	&+ \sum\limits_{j=0}^{\infty}\int_0^T U_j(\varrho(s))ds\left[\sum_{i=0}^{\infty}\int_0^T\partial_iV_j(\varrho(s))\mathcal{Z}_i(s,\theta^*) ds\right].
	\label{N-a12nm-a12-lim}
	\end{split}
\end{equation}
Similarly, as  $N,m,\frac{m}{\sqrt{N}}\rightarrow\infty$, $\sqrt{N}(b_2^{N,m}-b_2)$ converges in probability to
\begin{align}\label{N-b2nm-b2-lim}
		&
	\sum\limits_{j=0}^{\infty}\left[{\varrho}_j(T,\theta^*)-{\varrho}_j(0,\theta^*)\right]\Bigg[\sum_{i=0}^{\infty}\int_0^T\partial_iV_j(\varrho(s))\mathcal{Z}_i(s,\theta^*)ds\Bigg]\nonumber\\
		&+\sum\limits_{j=0}^{\infty}\left[\mathcal{Z}_j(T,\theta^*)-\mathcal{Z}_j(0,\theta^*)\right]\int_0^T V_j(\varrho(s))ds,
\end{align}
and $\sqrt{N}(a^{N,m}_{22}-a_{22})$ converges converges in probability to
\begin{equation}
	\begin{split}
	2\int_0^T V_j(\varrho(s))ds\Bigg[\sum\limits_{i=0}^{\infty}\int_0^T \partial_iV_j(s,\theta^*)\mathcal{Z}_i(s,\theta^*)ds\Bigg].
	\label{N-a22nm-a22}
	\end{split}
\end{equation}
Therefore, using $(\ref{N-b1nm-b1-lim})$--$(\ref{N-a22nm-a22})$, we deduce that
\begin{align*}
	&\mathcal{I}^{N,m}\xrightarrow{\text{p}}\mathcal{I},
\end{align*}	
Similarly, we obtain
\begin{align*}
	&\mathcal{J}^{N,m}\xrightarrow{\text{p}}\mathcal{J},\quad \mathcal{K}^{N,m}\xrightarrow{\text{p}}\mathcal{K}.
\end{align*}	

Finally, by Theorem \ref{LLN-Budhi}, we deduce that  $\mathcal{H}^{N,m}\rightarrow a_{11}a_{22}-(a_{12})^2$ in probability as $N,m\rightarrow\infty$.  Therefore, the proof is complete by the continuous mapping theorem (cf. \cite[Theorem 2.3]{van}).

\bibliographystyle{plain} 
\bibliography{biblio}

\begin{thebibliography}{10}

\bibitem{AMORINO2023}
C.~Amorino, A.~Heidari, V.~Pilipauskaite, and M.~Podolskij.
\newblock Parameter estimation of discretely observed interacting particle
  systems.
\newblock {\em Stochastic Processes and their Applications}, 163:350--386,
  2023.

\bibitem{Asanj+al2021}
A.~Asanjarani, Y.~Nazarathy, and P.~Taylor.
\newblock A survey of parameter and state estimation in queues.
\newblock {\em Queueing Syst}, 97:39--80, 2021.

\bibitem{belo2023}
D.~Belomestny, V.~Pilipauskaite, and M.~Podolskij.
\newblock Semiparametric estimation of mckean-vlasov sdes.
\newblock {\em Annales de l'Institut Henri Poincar\'e, Probabilit\'es et
  Statistiques}, 59(1):79--96, 2023.

\bibitem{Ben+Leboud2008}
M.~Bena\"{i}m and J.Y.~Le Boudec.
\newblock A class of mean field interaction models for computer and
  communication systems.
\newblock {\em Performance Evaluation}, 65(11--12):823--838, 2008.

\bibitem{Bishwal2011}
P.N. Bishwal.
\newblock Estimation in interacting diffusions: Continuous and discrete
  sampling.
\newblock {\em Applied Mathematics-a Journal of Chinese Universities Series B},
  02:1154--1158, 2011.

\bibitem{Budhi+Fried2019}
A.~Budhiraja and E.~Friedlander.
\newblock Diffusion approximations for load balancing mechanisms in cloud
  storage systems.
\newblock {\em Advances in Applied Probability}, 51(1):41--86, 2019.

\bibitem{Budh+Mukh+Wu2019}
A.~Budhiraja, D.~Mukherjee, and R.~Wu.
\newblock Supermarket model on graphs.
\newblock {\em The Annals of Applied Probability}, 29(3):1740--1777, 2019.

\bibitem{Curie1895}
P.~Curie.
\newblock Magnetic properties of materials at various temperatures.
\newblock {\em Ann. Chem. Phys}, 5(289), 1895.

\bibitem{Daw83}
D.A. Dawson.
\newblock Critical dynamics and fluctuations for a mean field model of
  cooperative behaviour.
\newblock {\em J. Statist. Phys.}, 41:29--85, 1983.

\bibitem{Daw2017}
D.A. Dawson.
\newblock Introductory lectures on stochastic population systems.
\newblock {\em arXiv:1705.03781 [math.PR]}, 2017.

\bibitem{Daw+Zhao2005}
D.A. Dawson, J.~Tang, and Y.Q. Zhao.
\newblock Balancing queues by mean field interactions.
\newblock {\em Queueing Syst}, 49:335--361, 2005.

\bibitem{Daw+Zhao2019}
D.A. Dawson, J.~Tang, and Y.Q. Zhao.
\newblock Performance analysis of joining the shortest queue model among a
  large number of queues.
\newblock {\em Asia-Pacific Journal of Operational Research}, 36(4), 2019.

\bibitem{Del+Four2016}
S.~Delattre and N.~Fournier.
\newblock Statistical inference versus mean field limit for hawkes processes.
\newblock {\em Electronic Journal of Statistics}, 10(1):1223--1295, 2016.

\bibitem{Dobru76}
R.L. Dobrushin and Y.M. Sukhov.
\newblock Asymptotic investigation of star-shaped message switching networks
  with a large number of radial rays.
\newblock {\em Probl. Inf. Trans.}, 12(1):49--66, 1976.

\bibitem{Eth+Kur86}
S.N. Ethier and T.G. Kurtz.
\newblock {\em The infinitely-many-alleles model with selection as a
  measure-valued diffusion. Stochastic Methods in Biology}, volume~70.
\newblock Springer-Verlag, Berlin-Heidelberg-N.Y., 1987.

\bibitem{Gart88}
J.~G\"artner.
\newblock On the mckean-vlasov limit for interacting diffusions.
\newblock {\em Math. Nachr.}, 137:197--248, 1988.

\bibitem{Genon+al2021}
V.~Genon-Catalot and C.~Lar{\'e}do.
\newblock Parametric inference for small variance and long time horizon
  mckean-vlasov diffusion models.
\newblock {\em Electronic Journal of Statistics}, 15(2):5811--5854.

\bibitem{Giesecke+al2019}
K.~Giesecke, G.~Schwenkler, and J.A. Sirignano.
\newblock Inference for large financial systems.
\newblock {\em Math. Fin.}, 30:823--838, 2019.

\bibitem{Giesecke+al2015}
K.~Giesecke, K.~Spiliopoulos, R.B. Sowers, and J.A. Sirignano.
\newblock Large portfolio asymptotics for loss from default.
\newblock {\em Mathematical Finance}, 25(1):77--114, 2015.

\bibitem{Graha2000}
C.~Graham.
\newblock Chaoticity on path space for a queueing network with selection of the
  shortest queue amongst several.
\newblock {\em J. Appl. Prob.}, 37(1):198--211, 2000.

\bibitem{Gra+Mel93}
C.~Graham and S.~M\'el\'eard.
\newblock Propagation of chaos for a fully connected loss network with
  alternate routing.
\newblock {\em Stoch. Processes Appl.}, 44:159--180, 1993.

\bibitem{Gra+Mel95}
C.~Graham and S.~M\'el\'eard.
\newblock Dynamic asymptotic results for a generalized star-shaped ioss
  network.
\newblock {\em Ann. Appl. Prob.}, 5, 1995.

\bibitem{McKean66}
H.P.~McKean Jr.
\newblock A class of markov processes associated with nonlinear parabolic
  equations.
\newblock {\em Proc. Natl. Acad. Sci. USA}, 56(6):1907--1911, 1966.

\bibitem{McKean66(2)}
H.P.~McKean Jr.
\newblock Speed of approach to equilibrium for kac's caricature of a maxwellian
  gas.
\newblock {\em Arch. Ration. Mech. Anal.}, 21(5):343--367, 1966.

\bibitem{Kac56}
M.~Kac.
\newblock Foundations of kinetic theory.
\newblock In Calif University~of California~Press, Berkeley, editor, {\em
  Proceedings of the Third Berkeley Symposium on Mathematical Statistics and
  Probability, Volume 3: Contributions to Astronomy and Physics}, pages
  171--197, 1956.

\bibitem{Kasonga90}
R.A. Kasonga.
\newblock Maximum likelihood theory for large interacting systems.
\newblock {\em SIAM Journal on Applied Mathematics}, 50(3):865--875, 1990.

\bibitem{Kley+Klu2015}
O.~Kley, C.~Kl\"{u}ppelberg, and L.~Reichel.
\newblock Systemic risk through contagion in a core-periphery structured
  banking network.
\newblock {\em Banach Center publications}, 104:133--149, 2015.

\bibitem{Liu2020}
C.~Liu.
\newblock Statistical inference for a partially observed interacting system of
  hawkes processes.
\newblock {\em Stochastic Processes and their Applications}, 130(9):5636--5694,
  2020.

\bibitem{Dell+Hoff2022}
L.~Della Maestra and M.~Hoffmann.
\newblock Nonparametric estimation for interacting particle systems:
  Mckean-vlasov models.
\newblock {\em Probab. Theory Relat. Fields}, 182:551--613, 2022.

\bibitem{Mel+Ben2015}
S.~M\'el\'eard and V.~Bansaye.
\newblock {\em Some Stochastic Models for Structured Populations: Scaling
  Limits and Long Time Behavior}.
\newblock Springer, 2015.

\bibitem{Mitz96}
M.~Mitzenmacher.
\newblock The power of two choices in randomized load balancing.
\newblock {\em PhD thesis, University of California at Berkeley}, 1996.

\bibitem{SHARROCK2023}
L.~Sharrock, N.~Kantas, P.~Parpas, and G.A. Pavliotis.
\newblock Online parameter estimation for the mckean-vlasov stochastic
  differential equation.
\newblock {\em Stochastic Processes and their Applications}, 162:481--546,
  2023.

\bibitem{Steele2004}
J.M. Steele.
\newblock {\em The Cauchy-Schwarz Master Class}.
\newblock Cambridge University Press, 2004.

\bibitem{van}
A.~W. van~der Vaart.
\newblock {\em Asymptotic Statistics}.
\newblock Cambridge Series in Statistical and Probabilistic Mathematics.
  Cambridge University Press, 1998.

\bibitem{Dobru96}
N.D. Vvedenskaya, R.L. Dobrushin, and F.I. Karpelevich.
\newblock Queueing system with selection of the shortest of two queues: An
  asymptotic approach.
\newblock {\em Problems Inform. Transmission}, 32(1):15--27, 1996.

\bibitem{Weiss1907}
P.~Weiss.
\newblock L'hypothese du champ mol\'eculaire et la propri\'et\'e
  ferromagnetique.
\newblock {\em J. Phys. Theor. Appl.}, 6(1):661--690, 1907.

\bibitem{Zhao2022}
Y.Q. Zhao.
\newblock Statistical inference for mean-field queueing models.
\newblock {\em Queueing Syst}, 100:569--571, 2022.

\end{thebibliography}
\end{document}